\documentclass[preprint,3p]{elsarticle}
\usepackage{cite}
\usepackage{amsmath}
\usepackage{graphicx,wrapfig}
\usepackage{amssymb,amsgen,color}
\usepackage{amsfonts}
\usepackage{hyperref}
\usepackage{multirow}
\usepackage{epsfig,epstopdf,color,bm}

\usepackage{braket,amsfonts}
\usepackage{booktabs}
\newcommand{\be}{\begin{equation}}
\newcommand{\ee}{\end{equation}}
\newcommand{\nn}{\nonumber}

\newtheorem{theorem}{Theorem}[section]
\newtheorem{lemma}{Lemma}[section]
\newtheorem{remark}{Remark}[section]

\biboptions{numbers,sort&compress}
\makeatletter 
\@addtoreset{equation}{section}
\@addtoreset{figure}{section}
\@addtoreset{table}{section}
\makeatother  

\graphicspath{{figures/}}

\begin{document}
	\begin{frontmatter}
		
\title{Uniform error bounds of exponential wave integrator methods for the long-time dynamics of the Dirac equation with small potentials}

\author[mymainaddress]{Yue Feng}
\ead{fengyue@u.nus.edu}

\author[mymainaddress,mysecondaryaddress]{Jia Yin\corref{mycorrespondingauthor}}
\cortext[mycorrespondingauthor]{Corresponding author}
\ead{jiayin@lbl.gov}

\address[mymainaddress]{Department of Mathematics, National University of Singapore, Singapore 119076, Singapore}
\address[mysecondaryaddress]{Computational Research Division, Lawrence Berkeley National Laboratory, Berkeley,\\ CA 94720, USA}
\begin{abstract}
Two exponential wave integrator Fourier pseudospectral (EWI-FP) methods are presented and analyzed for the long-time dynamics of the Dirac equation with small potentials characterized by $\varepsilon \in (0, 1]$ a dimensionless parameter. Based on the (symmetric) exponential wave integrator for temporal derivatives in phase space followed by applying the Fourier pseudospectral discretization for spatial derivatives, the EWI-FP methods are explicit and of spectral accuracy in space and second-order accuracy in time for any fixed $\varepsilon = \varepsilon_0$. Uniform error bounds are rigorously carried out at $O(h^{m_0}+\tau^2)$ up to the time at $O(1/\varepsilon)$ with the mesh size $h$, time step $\tau$ and $m_0$ an integer depending on the regularity of the solution. Extensive numerical results are reported to confirm our error bounds and comparisons of two methods are shown. Finally, dynamics of the Dirac equation in 2D are presented to validate the numerical schemes.
\end{abstract}


\end{frontmatter}

\section{Introduction}
In this paper, we consider the Dirac equation on the unit torus in one or two dimensions (1D or 2D), which can be represented in the two-component form with wave function $\Phi : = \Phi(t, {\mathbf{x}}) = (\phi_1(t, {\mathbf{x}}), \phi_2(t, {\mathbf{x}}))^T \in \mathbb{C}^2$ \citep{BCJT,CW,Dirac1,Dirac2}
\begin{equation}
\label{eq:Dirac_21}
i\partial_t\Phi =  \left(- i\sum_{j = 1}^{d}
	\sigma_j\partial_j + \sigma_3 \right)\Phi+ \varepsilon \left(V(t, \mathbf{x})I_2 - \sum_{j = 1}^{d}A_j(t, \mathbf{x})\sigma_j\right)\Phi, 
\end{equation}
where $i = \sqrt{-1}$, $t \ge 0$ is time, $\mathbf{x} = (x_1, \cdots, x_d)^T \in \mathbb{T}^d (d = 1, 2)$, $\partial_j = \frac{\partial}{\partial x_j} (j = 1, \cdots, d)$, $\varepsilon \in (0, 1]$ is a dimensionless parameter, $V:=V(t, \mathbf{x}) \in \mathbb{R}$ is the electric potential and $\mathbf{A}:=\mathbf{A}(t, \mathbf{x}) = (A_1(t, \mathbf{x}), A_2(t, \mathbf{x}), A_3(t, \mathbf{x}))^T \in \mathbb{R}^3$ stands for the magnetic potential. $I_2$ is the $2 \times 2$ identity matrix, and $\sigma_1, \sigma_2, \sigma_3$ are the Pauli matrices defined as
\begin{equation}
\label{Pauli}
\sigma_1 = \begin{pmatrix} 0 &\ \  1\\ 1 &\ \  0 \end{pmatrix}, \quad
\sigma_2 = \begin{pmatrix} 0 & \ \ -i \\ i &\ \   0\end{pmatrix}, \quad
\sigma_3 = \begin{pmatrix} 1 &\ \ 0 \\ 0 &\ \  -1 \end{pmatrix}.
\end{equation}
In order to study the dynamics of the Dirac equation \eqref{eq:Dirac_21}, the initial condition is taken as 
\begin{equation}
\label{eq:initial}
\Phi(t=0, \mathbf{x}) = \Phi_0(\mathbf{x}), \quad {\bf x} \in \mathbb{T}^d,\quad d = 1, 2.
\end{equation} 

For the Dirac equation \eqref{eq:Dirac_21} with $\varepsilon = 1$, i.e., the classical regime, there have been comprehensive analytical and numerical results in the literatures. Along the analytical front, for the existence and multiplicity of bound states and/or standing wave solutions, we refer to \citep{Das1,Das2,ES,GGT,Gross} and references therein. In the numerical aspect, different numerical methods have been proposed and analyzed, such as the finite difference time domain (FDTD) methods \citep{BCJT, BHM,MY}, exponential wave integrator Fourier pseudospectral (EWI-FP) method \citep{BCJT,BCJY}, time-splitting Fourier pseudospectral (TSFP) method \citep{BCY,FLB} and Gaussian beam method \citep{WHJY}. For more details related to the numerical schemes, we refer to \citep{AH,AL,BSG,CW,FLB,FLB1,Gosse,GSX} and references therein. 

On the other hand, when $0 < \varepsilon \leq 1$, one problem is to study the long-time dynamics of the Dirac equation \eqref{eq:Dirac_21} for $t \in [0, T_{\varepsilon}]$ with $T_{\varepsilon} = O(1/\varepsilon)$, i.e., adapt different numerical methods to get the long-time simulations and analyze how the errors of numerical schemes depend on the parameter $\varepsilon$ for $t \in [0, T_{\varepsilon}]$.  The other problem is to investigate  for the fixed $\varepsilon \in (0, 1]$, how the errors of different numerical schemes perform for $t \in [0, T_0]$ with large $T_0$. To the best of our knowledge, the long-time dynamics and conservation properties of the numerical schemes for the Schr\"odinger equation with small potential and nonlinear Schr\"odinger equation with small initial data have been studied in the literature \citep{CG,DF1,DF2,GL}. However, there are few results on error bounds of numerical methods for the long-time dynamics of the Dirac equation \eqref{eq:Dirac_21}. In our recent work, we adapt the finite difference method to discretize the Dirac equation \eqref{eq:Dirac_21} in time combined with different spatial discretizations. Long-time error bounds of the finite difference time domain (FDTD) methods and finite difference Fourier pseudospectral (FDFP) methods have been rigorously established for the Dirac equation \eqref{eq:Dirac_21} up to the time at $O(1/\varepsilon)$  with particular attentions paid to how the error bounds depend explicitly on the mesh size $h$ and the time step $\tau$ as well as the small parameter $\varepsilon \in (0, 1]$. Based on our error estimates, in order to obtain ``correct'' numerical approximations of the Dirac equation \eqref{eq:Dirac_21} up to the time at $O(1/\varepsilon)$, the $\varepsilon$-scalability (or meshing strategy requirement) of FDTD methods should be taken as
\begin{equation*}
h = O(\varepsilon^{1/2})\quad \mbox{and}\quad \tau = O(\varepsilon^{1/2}), \quad 0 < \varepsilon \leq 1.
\end{equation*}
Comparatively, the spatial error bounds of FDFP methods are uniform for any $\varepsilon \in (0, 1]$, which indicate that FDFP methods have better spatial resolution than FDTD methods to solve the Dirac equation with small potentials in the long-time regime. However, the temporal resolution of FDTD and FDFP methods depends on the small parameter $\varepsilon$, which causes severe numerical burdens as $\varepsilon \to 0^+$.

As we know, the exponential wave integrator Fourier pseudospectral (EWI-FP) method has been widely used to numerically solve dispersive partial differential equations (PDEs) in different regimes \citep{BC,CCO,GH,HLW,HL,HLS,HO}. The aim of this paper is to establish uniform error bounds of the EWI methods for the long-time dynamics of the Dirac equation \eqref{eq:Dirac_21} up to the time at $O(1/\varepsilon)$. We adapt the EWI-FP/symmetric EWI-FP (sEWI-FP) method to solve the Dirac equation with small potentials in the long-time regime and  uniform error bounds at $O(h^{m_0} + \tau^2)$ with $h$ mesh size, $\tau$ time step and $m_0$ depending on the regularity of the exact solution are proved. It is clear that the error bounds are independent of $\varepsilon$, which suggests that the numerical schemes are uniformly accurate for $\varepsilon \in (0, 1]$. The error estimates show that the $\varepsilon$-scalability of EWI-FP/sEWI-FP  method should be taken as $h = O(1)$ and $\tau = O(1)$ for the long-time dynamics of the Dirac equation with small potentials. Thus, EWI-FP/sEWI-FP  method offers compelling advantages over FDFP methods in temporal resolution when $ 0 < \varepsilon \ll 1$. 

The rest of this paper is organized as follows. In Section 2, the EWI-FP and sEWI-FP methods with some properties for the long-time dynamics of the Dirac equation with small potentials are presented. Error estimates for EWI-FP/sEWI-FP method are shown in Section 3. Extensive numerical examples in 1D and 2D are reported in Section 4. Finally, some conclusions are drawn in Section 5. Throughout this paper, we adopt the notation $A \lesssim B$ to represent that there exists a generic constant $C > 0$, which is independent of the mesh size $h$, time step $\tau$ and $\varepsilon$ such that $|A| \leq C B$.

\section{Numerical methods} 
In the section, we present the EWI-FP and sEWI-FP methods to solve the Dirac equation \eqref{eq:Dirac_21}. For simplicity of notations, we only present the numerical methods and their analysis in 1D. Generalization to higher dimensions is straightforward and the results are valid without modifications. In 1D, the Dirac equation \eqref{eq:Dirac_21} on the computation domain $\Omega = (a, b)$ with periodic boundary conditions collapses to 
\begin{align}
\label{eq:Dirac_1D}
&i\partial_t\Phi =  (- i \sigma_1 \partial_x + \sigma_3 )\Phi+ \varepsilon(V(t, x)I_2 - A_1(t, x)\sigma_1)\Phi, \ x \in \Omega,\ t > 0,\\
\label{eq:ib}
&\Phi(t, a) = \Phi(t, b),\ \partial_x \Phi(t, a) = \partial_x \Phi(t, b),\ t \geq 0; \ \Phi(0, x) = \Phi_0(x),\ x \in \overline{\Omega}, 
\end{align}
where $\Phi := \Phi(t, x)$, $\Phi_0(a) = \Phi_0(b)$ and $\Phi_0'(a) = \Phi_0'(b)$.

\subsection{The EWI-FP method}
Choose the time step $\tau = \Delta t > 0$ and mesh size $h := \Delta x = (b-a)/N$ with $N$ being an even positive integer. Denote time steps by $t_n = n\tau$ for $n \geq 0$ and grid points by $x_j := a+j h$, $j = 0, 1, \cdots, N$. Denote $X_N : = \{ U = (U_0, U_1, \cdots, U_N)^T \ | \ U_j \in \mathbb{C}^2,\  j =0, 1, \cdots, N,\  U_0 = U_N \}$ with the $l^2$-norm 
\begin{equation*}
\|U\|^2_{l^2} = h \sum^{N-1}_{j=0} |U_j|^2, \quad U \in X_N.
\end{equation*}
Define the index set $\mathcal{T}_N = \{l | l = -N/2,-N/2+1, \cdots, N/2-1\}$ and 
\begin{equation*}
Y_N = Z_N \times Z_N, \quad Z_N = \mbox{span}\{\phi_l(x) = e^{i\mu_l(x-a)},\ \mu_l = \frac{2 \pi l}{b-a},\  l \in \mathcal{T}_N\}.
\end{equation*}
Let $(C_p(\overline{\Omega}))^2$ be the function space consisting of all continuous periodic vector functions $U(x): \overline{\Omega} = [a, b] \to \mathbb{C}^2$. For any $U(x) \in (C_p(\overline{\Omega}))^2$ and $U \in X_N$, define $P_N: (L^2(\Omega))^2 \to Y_N$ as the standard projection operator \citep{ST}, $I_N : (C_p(\overline{\Omega}))^2 \to Y_N$ and $I_N : X_N \to Y_N$ as the standard interpolation operator, i.e.,
\begin{equation}
	(P_N U)(x) = \sum_{l \in \mathcal{T_N}}\widehat{U}_l e^{i\mu_l(x-a)}, \ (I_N U)(x) = \sum_{l \in \mathcal{T_N}}\widetilde{U}_l e^{i\mu_l(x-a)},\ a \leq x \leq b,
\end{equation} 
with 
\begin{equation}
\label{eq:Fc}
\widehat{U}_l = \frac{1}{b-a}\int^b_a U(x) e^{-i\mu_l(x-a)} dx,\quad \widetilde{U}_l = \frac{1}{N}\sum_{j=0}^{N-1}U_j e^{-2ijl\pi/N},\quad l \in \mathcal{T}_N, 	
\end{equation}
where $U_j = U(x_j)$ when $U$ is a function.

The Fourier spectral discretization for the Dirac equation \eqref{eq:Dirac_1D} is to find $\Phi_N(t, x) \in Y_N$, i.e.,
\begin{equation}
\label{eq:Phi_F}
\Phi_N(t, x) = \sum_{l \in \mathcal{T}_N} \widehat{(\Phi_N)}_l(t)e^{i\mu_l(x-a)}, \quad  a \leq x \leq b,\quad t\geq 0,
\end{equation}
such that, for $a < x <b$ and $t > 0$, $\Phi_N :=\Phi_N(t, x)$, satisfies
\begin{equation}
i \partial_t\Phi_N = \left( -i \sigma_1 \partial_x + \sigma_3\right)\Phi_N +\varepsilon P_N(V\Phi_N) -\varepsilon \sigma_1P_N(A_1\Phi_N),
\label{eq:Phi_eq}
\end{equation}
where we take $V := V(t, x)$, $A_1 := A_1(t, x)$ in short.
Plugging \eqref{eq:Phi_F} into \eqref{eq:Phi_eq}, by noticing the orthogonality of the Fourier basis functions, we derive 
\begin{equation*}
i\frac{d}{dt}\widehat{(\Phi_N)}_l(t) = \left(\mu_l \sigma_1 + \sigma_3\right)\widehat{(\Phi_N)}_l(t)+ \varepsilon\widehat{(V\Phi_N)}(t) - \varepsilon \sigma_1 \widehat{(A_1\Phi_N)}(t).
\end{equation*}
For each $l \in \mathcal{T}_N$, when $t$ is near $t = t_n$ $(n \geq 0)$, we rewrite the above ODEs as 
\begin{equation}
i\frac{d}{ds}\widehat{(\Phi_N)}_l(t_n + s) =\Gamma_l \widehat{(\Phi_N)}_l(t_n + s)+ \varepsilon \widehat{F}^n_l(s), \quad s \in \mathbb{R},
\label{eq:ODEs}
\end{equation}
where
\begin{equation}
\widehat{F}^n_l(s) = \widehat{(G\Phi_N)_l}(t_n+s),\quad G(t, x) = V(t, x)I_2 -\sigma_1 A_1(t, x),\quad s,t \in \mathbb{R},
\label{eq:FG}
\end{equation}
and $\Gamma_l = \mu_l \sigma_1+\sigma_3 = Q_l D_l (Q_l)^{\ast}$  with $\delta_l =\sqrt{1+\mu^2_l}$ and 
\begin{equation}
\Gamma_l = \begin{pmatrix} 1 &\  \mu_l \\ \mu_l &\  -1 \end{pmatrix}, \
 Q_l = \begin{pmatrix} \frac{1+\delta_l}{\sqrt{2\delta_l(1+\delta_l)}}  & \ -\frac{\mu_l}{\sqrt{2\delta_l(1+\delta_l)}} \\ \frac{\mu_l}{\sqrt{2\delta_l(1+\delta_l)}} &\    \frac{1+\delta_l}{\sqrt{2\delta_l(1+\delta_l)}} \end{pmatrix}, \
 D_l = \begin{pmatrix} \delta_l &\ 0 \\ 0 &\   -\delta_l \end{pmatrix}.
\label{eq:GammaQD}
\end{equation}
Solving the above ODE \eqref{eq:ODEs} via the integrating factor method, we obtain
\begin{equation}
\widehat{(\Phi_N)}_l(t_n + s) = e^{-is\Gamma_l}\widehat{(\Phi_N)}_l(t_n) - i\varepsilon \int_0^s e^{i(w-s)\Gamma_l}\widehat{F}^n_l(w) d w,\quad s\in \mathbb{R}.
\label{eq:ifm}
\end{equation}
Taking $s = \tau$ in \eqref{eq:ifm}, we arrive at
\begin{equation}
\widehat{(\Phi_N)}_l(t_{n+1}) = 	e^{-i\tau\Gamma_l}\widehat{(\Phi_N)}_l(t_n) - i\varepsilon \int_0^{\tau} e^{i(w-\tau)\Gamma_l}\widehat{F}^n_l(w) dw.
\label{eq:phi_n1}
\end{equation}
To obtain an explicit numerical scheme with second order accuracy in time, we approximate the integral in \eqref{eq:phi_n1} via the Gautschi-type rule which has been widely used for integrating the ODEs \citep{Gau,Gr1,Gr2,HO} as following
\begin{equation}
\int_0^{\tau} e^{i(w-\tau)\Gamma_l}\widehat{F}^0_l(w) d w  \approx \int_0^{\tau} e^{i(w-\tau)\Gamma_l}d w \widehat{F}^0_l(0)  = -i \Gamma_l^{-1}\left[I_2-e^{-i\tau\Gamma_l}\right] \widehat{F}^0_l(0),
\label{eq:approx1}
\end{equation}
and for $n \geq 1$,
\begin{align}
\int_0^{\tau} e^{i(w-\tau)\Gamma_l}\widehat{F}^n_l(w) d w	 \approx & \int_0^{\tau} e^{i(w-\tau)\Gamma_l}\left(\widehat{F}^n_l(0) + w \delta^-_t  \widehat{F}^n_l(0)  \right) d w \nonumber \\
 = & -i\Gamma_l^{-1}\left[I_2-e^{-i\tau\Gamma_l}\right]\widehat{F}^n_l(0) + \left[-i\tau \Gamma_l^{-1}+ \Gamma_l^{-2}(I_2-e^{-i\tau\Gamma_l})\right]\delta^-_t \widehat{F}^n_l(0),
\label{eq:approx2}
\end{align}
where we have approximated the time derivative $\partial_t \widehat{F}^n_l(s)$ at $s = 0$ by finite difference as
\begin{equation*}
\partial_t \widehat{F}^n_l(0) \approx \delta^-_t  \widehat{F}^n_l(0)  = \frac{1}{\tau}\left[ \widehat{F}^n_l(0) -  \widehat{F}^{n-1}_l(0) \right].
\end{equation*}
Then, we are going to describe our numerical scheme. Let $\Phi^n_N(x)$ be the approximation of $\Phi_N(t_n, x)$ for $n\geq 0$. Choose $\Phi^0_N(x)=(P_N \Phi_0)(x)$, then the exponential wave integrator Fourier spectral (EWI-FS) discretization for the Dirac equation \eqref{eq:Dirac_1D} is to update the numerical approximation $\Phi^{n+1}_N \in Y_N$ $(n =0 ,1, \cdots)$ as
\begin{equation}
\Phi^{n+1}_N(x) = \sum_{l \in \mathcal{T}_N} \widehat{(\Phi^{n+1}_N)_l}e^{i\mu_l(x-a)}, \quad a \leq x \leq b, \quad n \geq 0,	
\label{eq:Phi_hat}
\end{equation}
where for $l \in \mathcal{T}_N$,
\begin{equation}
\widehat{(\Phi^{n+1}_N)_l} = \left\{
\begin{split}
& e^{-i\tau\Gamma_l}\widehat{(\Phi^{0}_N)}_l - \varepsilon \Gamma_l^{-1}\left[I_2-e^{-i\tau\Gamma_l}\right] \widehat{(G(t_0)\Phi^0_N)}_l, \quad & n = 0, \\
&e^{-i\tau\Gamma_l}\widehat{(\Phi^{n}_N)}_l - i\varepsilon Q^{(1)}_l(\tau)\widehat{(G(t_n)\Phi^n_N)}_l - i\varepsilon Q^{(2)}_l(\tau)\delta^-_t\widehat{(G(t_n)\Phi^n_N)}_l, \quad & n \geq 1,
\end{split}\right.
\label{eq:Phi_hat_co}
\end{equation}
with the matrices $Q^{(1)}_l(\tau)$ and $Q^{(2)}_l(\tau)$ given as
\begin{equation*}
Q^{(1)}_l(\tau) = -i \Gamma_l^{-1}\left[I_2-e^{-i\tau\Gamma_l}\right],\ Q^{(2)}_l(\tau) = -i\tau \Gamma^{-1}_l + \Gamma^{-2}_l\left[I_2-e^{-i\tau\Gamma_l}\right].
\end{equation*}

The above procedure is not suitable in practice due to the difficulty of computing the integrals in \eqref{eq:Fc}. We now present an efficient implementation by choosing $\Phi^0_N(x)$ as the interpolation of $\Phi_0(x)$ on the grids $\{x_j, \ j=0, 1, \cdots, N\}$, i.e., $\Phi^0_N(x) = (I_N\Phi_0)(x)$ and approximate the integrals by a quadrature rule on the grids.

Let $\Phi^n_j$ be the numerical approximation of $\Phi(t_n, x_j)$ for $j = 0, 1, \cdots, N$ and $n \geq 0$, and denote $\Phi^n \in X_N$ as the vector with components $\Phi^n_j$. Choosing $\Phi^0_j = \Phi_0(x_j)$ $ (j = 0, 1, \cdots, N)$, the exponential wave integrator Fourier pseudospectral (EWI-FP) method for computing $\Phi^{n+1}$ for $n \geq 0$ reads
\begin{equation}
\Phi^{n+1}_j(x) = \sum_{l \in \mathcal{T}_N} \widetilde{(\Phi^{n+1})_l}e^{2ijl\pi/N}, \quad j= 0 ,1, \cdots, N,	
\label{eq:Phi_tilde}
\end{equation}
where 
\begin{equation}
\widetilde{(\Phi^{n+1})_l} = \left\{
\begin{split}
& e^{-i\tau\Gamma_l}\widetilde{(\Phi^{0})}_l - \varepsilon \Gamma_l^{-1}\left[I_2-e^{-i\tau\Gamma_l}\right] \widetilde{(G(t_0)\Phi^0)}_l, \ & n = 0, \\
&e^{-i\tau\Gamma_l}\widetilde{(\Phi^{n})}_l - i\varepsilon Q^{(1)}_l(\tau)\widetilde{(G(t_n)\Phi^n)}_l  - i\varepsilon Q^{(2)}_l(\tau)\delta^-_t\widetilde{(G(t_n)\Phi^n)}_l, \ & n \geq 1.
\end{split}\right.
\label{eq:Phi_tilde_co}
\end{equation}

\subsection{The sEWI-FP method}
We note that $e^{is\Gamma_l} = \cos(s\Gamma_l) + i \sin(s \Gamma_l)$ $(s\in \mathbb{R}, l \in \mathcal{T}_N)$ and
\begin{equation*}
\sin(s\Gamma_l) = Q_l \begin{pmatrix} \sin(s \delta_l) &\ \  0 \\ 0 &\ \  -\sin(s \delta_l) \end{pmatrix} Q^{\ast}_l, \quad \cos(s \Gamma_l) = \cos(s \delta_l) I_2.	
\end{equation*}
For $n \geq 1$, taking $s = \tau$ and $s = -\tau$ in \eqref{eq:ifm}, and subtracting one from the other, we can obtain 
\begin{align}
\widehat{(\Phi_N)}_l(t_{n+1}) = & \  \widehat{(\Phi_N)}_l(t_{n-1}) - 2i\sin(\tau \Gamma_l)\widehat{(\Phi_N)}_l(t_n)  -  i\varepsilon \int_0^{\tau} \cos((w-\tau)\delta_l)\left(\widehat{F}^n_l(w) + \widehat{F}^n_l(-w)\right) d w \nonumber\\
& \ + \varepsilon \int_0^{\tau} \sin((w-\tau)\Gamma_l)\left(\widehat{F}^n_l(w) - \widehat{F}^n_l(-w)\right) d w.
\label{eq:sEWI2}
\end{align}
Similarly, we approximate the integrals in \eqref{eq:sEWI2} via the Gautschi-type/trapezoidal rules. For $n = 0$, we use the same approximation as \eqref{eq:approx1} and for $ n \geq 1$,
\begin{align*}
& \int^{\tau}_0 \cos\left((w-\tau)\delta_l\right) \left(\widehat{F}^n_l(w) + \widehat{F}^n_l(-w) \right) d w \\
\approx & \int^{\tau}_0 \cos\left((w-\tau)\delta_l\right) \left(2\widehat{F}^n_l(0) + (w-w) \partial_t \widehat{F}^n_l(0) \right) d w =  \frac{2}{\delta_l}\sin(\tau \delta_l) \widehat{F}^n_l(0), \\
& \int^{\tau}_0 \sin\left((w-\tau)\Gamma_l\right) \left(\widehat{F}^n_l(w) - \widehat{F}^n_l(-w) \right) d w \\
\approx &\  \frac{\tau}{2}\left[\sin(-\tau\Gamma_l) {\bf{0}} + \sin(0 \Gamma_l)\left(\widehat{F}^n_l(\tau) - \widehat{F}^n_l(-\tau) \right)\right] = \bf{0}.
\end{align*}
Then, by choosing $\Phi^0_N(x) = (P_N \Phi_0)(x)$, a symmetric exponential wave integrator Fourier spectral (sEWI-FS) discretization for the Dirac equation \eqref{eq:Dirac_1D} is to update the numerical approximation $\Phi^{n+1}_N \in Y_N$ $(n =0 ,1, \cdots)$ as
\begin{equation}
\Phi^{n+1}_N(x) = \sum_{l \in \mathcal{T}_N} \widehat{(\Phi^{n+1}_N)_l}e^{i\mu_l(x-a)}, \quad a \leq x \leq b, \quad n \geq 0,	
\label{eq:sPhi_hat}
\end{equation}
where for $l \in \mathcal{T}_N$,
\begin{equation}
\widehat{(\Phi^{n+1}_N)_l} = \left\{
\begin{split}
& e^{-i\tau\Gamma_l}\widehat{(\Phi^{0}_N)}_l - \varepsilon \Gamma_l^{-1}\left[I_2-e^{-i\tau\Gamma_l}\right] \widehat{(G(t_0)\Phi^0_N)}_l, & n = 0, \\
& \widehat{(\Phi^{n-1}_N)}_l -2i\sin(\tau\Gamma_l) \widehat{(\Phi^n_N)}_l - i\frac{2\varepsilon}{\delta_l} \sin(\tau\delta_l) \widehat{(G(t_n)\Phi^n_N)}_l, & n \geq 1.
\end{split}\right.
\label{eq:sPhi_hat_co}
\end{equation}
For $n \geq 1$, this scheme is unchanged if we interchange $n + 1 \leftrightarrow n -1$ and $\tau \leftrightarrow - \tau$.
Similarly, a symmetric exponential wave integrator Fourier pseudospectral (sEWI-FP) discretization for the Dirac equation \eqref{eq:Dirac_1D} is computing $\Phi^{n+1}$ $(n =0 ,1, \cdots)$ as
\begin{equation}
\Phi^{n+1}_j = \sum_{l \in \mathcal{T}_N} \widetilde{(\Phi^{n+1})_l}e^{2ijl\pi/N}, j = 0, 1, \cdots, N,	
\label{eq:sEWI_Phi}
\end{equation}
where for $l \in \mathcal{T}_N$,
\begin{equation}
\widetilde{(\Phi^{n+1})_l} = \left\{
\begin{split}
& e^{-i\tau\Gamma_l}\widetilde{(\Phi_0)}_l - \varepsilon \Gamma_l^{-1}\left[I_2-e^{-i\tau\Gamma_l}\right] \widetilde{(G(t_0)\Phi_0)}_l, & n = 0, \\
& \widetilde{(\Phi^{n-1})}_l  -2i\sin(\tau\Gamma_l) \widetilde{(\Phi^n)}_l - i\frac{2\varepsilon}{\delta_l} \sin(\tau\delta_l) \widetilde{(G(t_n)\Phi^n)}_l, & n \geq 1.
\end{split}\right.
\label{eq:sEWI_Phitilde}
\end{equation}

\begin{remark}
The EWI-FP \eqref{eq:Phi_tilde}-\eqref{eq:Phi_tilde_co} and sEWI-FP \eqref{eq:sEWI_Phi}-\eqref{eq:sEWI_Phitilde} are explicit and can be solved efficiently by the fast Fourier transform (FFT). The memory cost is $O(N)$ and the computational cost per time step is $O(N \ln N)$. 	
\end{remark}

\subsection{Stability conditions}
To consider the linear stability of the EWI-FP and sEWI-FP methods, we assume that $A_1(t, x) \equiv A^0_1$ and $V(t, x) \equiv V^0$ with $A^0_1$ and $V^0$ two real constants in the Dirac equation \eqref{eq:Dirac_1D}. In this case, we adopt the standard von Neumann stability analysis \citep{Smith} to prove that the errors grow exponentially at most. 
\begin{lemma}
\label{lemma:EWI_sta}
The EWI-FP \eqref{eq:Phi_tilde}-\eqref{eq:Phi_tilde_co} is stable under the stability condition
\begin{equation*}
0 < \tau \lesssim 1,
\end{equation*}
for any given $0 < \varepsilon \leq 1$.
\end{lemma}
\noindent
\emph{Proof.}
Plugging 
\begin{equation}
\label{eq:phi_c}
\Phi^{n}_j = \sum_{l \in \mathcal{T}_N} \xi^n_l \widetilde{(\Phi^0)}_l e^{i\mu_l(x_j-a)} = \sum_{l \in \mathcal{T}_N}  \xi^n_l \widetilde{(\Phi^0)}_l e^{2ijl\pi/N}, \quad 0 \leq j \leq N,
\end{equation}	
into \eqref{eq:phi_n1} with the integration approximated by \eqref{eq:Phi_tilde_co}, $\xi_l$ and $\widetilde{(\Phi^0)}_l$ the amplification factor and the Fourier coefficient at $n=0$ of the $l$th mode in phase space, respectively, we obtain for $l \in \mathcal{T}_N$
\begin{equation}
\xi^2_l \widetilde{(\Phi^0)}_l = \xi_l e^{-i\tau\Gamma_l} \widetilde{(\Phi^0)}_l - i\varepsilon \int^{\tau}_0 e^{i(w-\tau)\Gamma_l} (V^0 I_2 - A^0_1\sigma_1) \left(\xi_l + \frac{w}{\tau}(\xi_l-1)\right)\widetilde{(\Phi^0)}_l d w.
\label{eq:linear1}
\end{equation}
Denoting $C = |V^0| + |A^0_1|$, taking the $l^2$ norms of the vectors on both sides of \eqref{eq:linear1} and then dividing by the $l^2$ norms of $(\widetilde{\Phi^0})_l$, in view of the properties of $e^{-i\tau\Gamma_l}$, we get
\begin{equation*}
|\xi_l|^2 \leq \left(1+ \varepsilon C\tau+\frac{\varepsilon C}{2}\tau\right)|\xi_l| +\frac{\varepsilon C}{2}\tau,	
\end{equation*}
which implies
\begin{equation*}
\left(|\xi_l| - \frac{1+3\varepsilon C\tau/2}{2}\right)^2\leq \frac{1+5\varepsilon C\tau +9\varepsilon ^2C^2\tau^2/4}{4} \leq \left(\frac{1+5\varepsilon C\tau/2}{2}\right)^2.
\end{equation*}
Thus, we obtain for $l \in \mathcal{T}_N$
\begin{equation*}
|\xi_l| \leq 1 + 2\varepsilon C\tau,
\end{equation*}
which implies that the EWI-FP \eqref{eq:Phi_tilde}-\eqref{eq:Phi_tilde_co} is stable.
\hfill $\square$ \bigskip

\begin{lemma}
The sEWI-FP \eqref{eq:sEWI_Phi}-\eqref{eq:sEWI_Phitilde} is stable under the stability condition
\begin{equation}
\label{eq:sEWI_con}
0 < \tau < \min\Big\{\frac{h\pi}{3\sqrt{h^2+\pi^2}}, \frac{2-\sqrt{3}}{2(|V^0| + |A^0_1|)}\Big\}	
\end{equation}
for any given $0 < \varepsilon \leq 1$.
\end{lemma}

\noindent
\emph{Proof.}
Similarly to the proof of Lemma \ref{lemma:EWI_sta}, noticing \eqref{eq:sEWI_Phitilde}, we can get
\begin{equation}
\xi^2_l \widetilde{(\Phi^0)}_l = \widetilde{(\Phi^0)}_l - 2i\xi_l\sin(\tau\Gamma_l)\widetilde{(\Phi^0)}_l -2i \xi_l \varepsilon \delta^{-1}_l \sin(\tau \delta_l)\left(V^0I_2-A^0_1\sigma_1\right)\widetilde{(\Phi^0)}_l. 
\label{eq:sEWI_sta1}
\end{equation}
Multiplying both sides of \eqref{eq:sEWI_sta1} by $\overline{\xi_l} \widetilde{(\Phi^0)}_l^{\ast}$ and then taking the real part and dividing both sides by $\widetilde{(\Phi^0)}_l^{\ast}\widetilde{(\Phi^0)}_l$, in view of Hermitian matrices $\Gamma_l$ and $\sigma_1$, we get
\begin{equation}
|\xi_l|^2 \text{Re}(\xi_l) = \text{Re}(\overline{\xi_l}),	
\end{equation}
which implies $|\xi_l| = 1$ if $\text{Re}(\xi_l) \neq 0$. On the other hand, if $\text{Re}(\xi_l) = 0$, we can take $\xi_l = i c_l$ with $c_l \in \mathbb{R}$, and \eqref{eq:sEWI_sta1} leads to
\begin{equation}
-c^2_l \widetilde{(\Phi^0)}_l = \widetilde{(\Phi^0)}_l + 2c_l\sin(\tau\Gamma_l)\widetilde{(\Phi^0)}_l + 2 c_l \varepsilon \delta^{-1}_l \sin(\tau \delta_l)\left(V^0I_2-A^0_1\sigma_1\right)\widetilde{(\Phi^0)}_l. 
\label{eq:sEWI_sta2}	
\end{equation}
Denoting $C = |V^0| + |A^0_1|$, multiplying both sides of \eqref{eq:sEWI_sta2} by $\widetilde{(\Phi^0)}_l^{\ast}$ and then dividing both sides by $\widetilde{(\Phi^0)}_l^{\ast}\widetilde{(\Phi^0)}_l$, noticing $\sin(\tau \delta_l) \leq \frac{\sqrt{3}}{2}$ under the stability constraint and $|\delta_l^{-1} \sin(\tau\delta_l)| \leq \tau$, we obtain
\begin{equation}
c_l^2 + 1 \leq \sqrt{3}|c_l| + 2\varepsilon \tau C|c_l|.
\end{equation}
If $\tau < \frac{2-\sqrt{3}}{2C}$, there is no real number $c_l$ satisfying the above inequality. It follows that the sEWI-FP \eqref{eq:sEWI_Phi}-\eqref{eq:sEWI_Phitilde} is stable under the stability condition \eqref{eq:sEWI_con}.
\hfill $\square$ \bigskip

\section{Uniform error bounds for long-time dynamics}
In this section, we rigorously carry out uniform error bounds of the EWI methods for the Dirac equation \eqref{eq:Dirac_1D} up to the time at $O(1/\varepsilon)$.
\subsection{Main results}
In order to obtain the error bounds for the EWI-FS method \eqref{eq:Phi_hat}-\eqref{eq:Phi_hat_co} and sEWI-FS method \eqref{eq:sPhi_hat}-\eqref{eq:sPhi_hat_co}, motivated by the results in \citep{BCJT, CC}, we assume that there exists an integer $m_0 \geq 2$ such that the exact solution $\Phi(t, x)$ of the Dirac equation \eqref{eq:Dirac_1D} up to the time at $t = T_0/\varepsilon $ satisfies
\begin{align*}
\textrm{(A)} \quad &\|\Phi\|_{L^{\infty}([0, T_0/\varepsilon]; (H^{m_0}_p)^2)} \lesssim 1,	\quad
\|\partial_t \Phi\|_{L^{\infty}([0, T_0/\varepsilon]; (L^2)^2)} \lesssim 1,\quad \|\partial_{tt} \Phi\|_{L^{\infty}([0, T_0/\varepsilon]; (L^2)^2)} \lesssim 1,
\end{align*}
where $H^k_p(\Omega) = \{u \ | \ u \in H^k(\Omega), \partial^l_x u(a) =\partial^l_x u(b), l = 0, \cdots, k-1\}$. In addition, we assume electromagnetic potentials satisfy
\begin{equation*}
\textrm{(B)} \quad \|V\|_{W^{2, \infty}([0, T_0/\varepsilon]; L^{\infty})} + \|A_1\|_{W^{2, \infty}([0, T_0/\varepsilon]; L^{\infty})} \lesssim 1.
\end{equation*}
We can establish the following error estimates for the EWI-FS and sEWI-FS methods.

\begin{theorem}
Let $\Phi^n_N(x)$ be the approximation obtained from the EWI-FS \eqref{eq:Phi_hat}-\eqref{eq:Phi_hat_co}. Under the assumptions (A) and (B), there exists $\tau_0 > 0$ sufficiently small and independent of $\varepsilon$ such that, for any $0 < \varepsilon \leq 1$, when $0 <\tau \leq \tau_0$, we have 
\begin{equation}
\left\|\Phi(t_n, x) - \Phi^n_N(x)\right\|_{L^2} \leq C_{T_0} \left(h^{m_0} + \tau^2\right), \quad 0 \leq n \leq \frac{T_0/\varepsilon}{\tau},
\label{eq:EWI_error}
\end{equation}
where $C_{T_0}$ is a constant depending on $T_0$ and independent of $h, \tau$ and $\varepsilon$.
\label{thm:EWI-FS}
\end{theorem}

\begin{theorem}
Let $\Phi^n_N(x)$ be the approximation obtained from the sEWI-FS \eqref{eq:sPhi_hat}-\eqref{eq:sPhi_hat_co}. Under the assumptions (A) and (B), there exists $\tau_0 > 0$ sufficiently small and independent of $\varepsilon$ such that, for any $0 < \varepsilon \leq 1$,  under the stability condition \eqref{eq:sEWI_con}, when $0 <\tau \leq \tau_0$, we have the following error estimate
\begin{equation}
\left\|\Phi(t_n, x) - \Phi^n_N(x)\right\|_{L^2} \leq C_{T_0} \left(h^{m_0} + \tau^2\right), \quad 0 \leq n \leq \frac{T_0/\varepsilon}{\tau},
\label{eq:sEWI_error}
\end{equation}
where $C_{T_0}$ is a constant depending on $T_0$ and independent of $h, \tau$ and $\varepsilon$.
\label{thm:sEWI-FS}
\end{theorem}

\begin{remark}
The same error estimates in Theorem \ref{thm:EWI-FS} and Theorem \ref{thm:sEWI-FS} also hold for the EWI-FP \eqref{eq:Phi_tilde}-\eqref{eq:Phi_tilde_co} and the sEWI-FP \eqref{eq:sEWI_Phi}-\eqref{eq:sEWI_Phitilde}. The proofs are quite similar to the EWI-FS and sEWI-FS methods and we just show the details for the EWI-FS and sEWI-FS methods for brevity.
\end{remark}

\subsection{Proof for Theorem \ref{thm:EWI-FS}}
Define the error function $\textbf{e}^n(x) \in Y_N$ for $n \geq 0$ as
\begin{equation}
\textbf{e}^n(x) : = P_N \Phi(t_n, x) - \Phi^n_N(x) = \sum_{l \in \mathcal{T}_N} \widehat{\textbf{e}}^n_l e^{i \mu_l(x-a)},\quad a \leq x \leq b.
\end{equation}
By the triangular inequality and standard projection result, we get
\begin{align*}
\|\Phi(t_n, x) - \Phi^n_N(x)\|_{L^2}& \leq \|\Phi(t_n, x) - P_N \Phi(t_n, x)\|_{L^2} + \|\textbf{e}^n(x)\|_{L^2} \\
& \leq C_0h^{m_0} +  \|\textbf{e}^n(x)\|_{L^2}, \quad 0 \leq n \leq \frac{T_0/\varepsilon}{\tau},
\end{align*}
where $C_0$ is a constant independent of $h$, $\tau$ and $\varepsilon$. Then we only need to estimate $\|\textbf{e}^n(x)\|_{L^2}$ for $0 \le n \le \frac{T_0/\varepsilon}{\tau}$. It is easy to see that $\|{\textbf e}^0\|_{L^2} = 0$ since $\Phi_N^0(x) = (P_N\Phi_0)(x)$.

Define the local truncation error $\xi^n(x) = \sum_{l \in \mathcal{T}_N} \widehat{\xi}^n_l e^{i \mu_l(x-a)} \in Y_N$ of the EWI-FS \eqref{eq:Phi_hat}-\eqref{eq:Phi_hat_co} as
\begin{equation}
\begin{split}
\widehat{\xi}^0_l = &\ \widehat{(\Phi(\tau))}_l -  e^{-i\tau\Gamma_l}\widehat{(\Phi(0))}_l +  \varepsilon \Gamma_l^{-1}\left[I_2-e^{-i\tau\Gamma_l}\right] \widehat{(G(0)\Phi(0))}_l, \\
\widehat{\xi}^n_l = &\ \widehat{(\Phi(t_{n+1}))}_l - e^{-i\tau\Gamma_l}\widehat{(\Phi(t_n))}_l + i\varepsilon Q^{(1)}_l(\tau)\widehat{(G(t_n)\Phi(t_n))}_l  + i\varepsilon Q^{(2)}_l(\tau)\delta^-_t\widehat{(G(t_n)\Phi(t_n))}_l,\ n \geq 1,
\end{split}
\label{eq:local}
\end{equation}
where we write $\Phi(t)$ and $G(t)$ in short for $\Phi(t, x)$ and $G(t, x)$, respectively. 

Firstly, we estimate the local truncation error $\xi^n(x)$. Multiplying both sides of the Dirac equation \eqref{eq:Dirac_1D} by $e^{i\mu_l(x-a)}$ and integrating over the interval $(a, b)$, we can get the equations for $\widehat{(\Phi(t))}_l$, which are exactly the same as \eqref{eq:ODEs} with $\Phi_N$ replaced by $\Phi(t, x)$. Replacing $\Phi_N$ by $\Phi(t, x)$, we use the same notations $\widehat{F}^n_l(s)$ as in \eqref{eq:FG} and the time derivatives of $\widehat{F}^n_l(s)$ enjoy the same properties of time derivatives of $\Phi(t, x)$. Thus, the same representation \eqref{eq:phi_n1} holds for $\widehat{(\Phi(t_n))}_l$ with $n \geq 1$. From the derivation of the EWI-FS method, it is clear that the error $\xi^n(x)$ comes from the approximations for the integrals in \eqref{eq:approx1} and \eqref{eq:approx2}. Thus, we have
\begin{equation}
\widehat{\xi}^0_l  = -i\varepsilon \int^{\tau}_0 e^{i(s-\tau)\Gamma_l}\left(\widehat{F}^0_l(s) - \widehat{F}^0_l(0)\right) d s = 	-i\varepsilon\int^{\tau}_0 \int^s_0 e^{i(s-\tau)\Gamma_l}\partial_{s_1} \widehat{F}^0_l(s_1) d s_1 d s,
\label{eq:xi0}
\end{equation}
and for $n \geq 1$
\begin{align}
\widehat{\xi}^n_l &= -i\varepsilon \int^{\tau}_0 e^{i(s-\tau)\Gamma_l}\left(\int^s_0 \partial_{s_1} \widehat{F}^n_l(s_1) d s_1 - \frac{s}{\tau}\int^\tau_0 \partial_{\theta} \widehat{F}^{n-1}_l(\theta)d\theta \right) d s \nonumber \\
& = -i\varepsilon \int^{\tau}_0 e^{i(s-\tau)\Gamma_l}\left(\int^s_0\int_0^{s_1} \partial_{s_2s_2} \widehat{F}^n_l(s_2)ds_2 d s_1 + s\int_0^1\int^\tau_{\theta\tau} \partial_{\theta_1\theta_1} \widehat{F}^{n-1}_l(\theta_1)d\theta_1d\theta\right) d s\
\label{eq:xin}
\end{align}
Subtracting \eqref{eq:Phi_hat_co} from \eqref{eq:local}, we obtain the following error function
\begin{align}
& \widehat{\mathbf{e}}^{n+1}_l = e^{-i\tau\Gamma_l} \widehat{\mathbf{e}}^{n}_l + \widehat{R}^n_l + \widehat{\xi}^n_l, \quad 1 \leq n \leq \frac{T_0/\varepsilon}{\tau} -  1,  \label{eq:errorfun} \\
& \widehat{\mathbf{e}}^0_l = \textbf{0}, \quad \widehat{\mathbf{e}}^1_l = \widehat{\xi}^0_l, \quad l \in \mathcal{T}_N, 	
\end{align}
where $R^n(x) = \sum_{l \in \mathcal{T}_N} \widehat{R}^n_l e^{i \mu_l(x-a)} \in Y_N$ for $ n \geq 1$, with $\widehat{R}^n_l$ given by 
\begin{equation}
\widehat{R}^n_l = \  -i\varepsilon Q^{(1)}_l(\tau)\left[\widehat{(G(t_n)\Phi(t_n))}_l - \widehat{(G(t_n)\Phi^n_N)}_l\right] -i \varepsilon Q^{(2)}_l(\tau)\left[\delta^-_t\widehat{(G(t_n)\Phi(t_n))}_l - \delta^-_t \widehat{(G(t_n)\Phi^n_N)}_l\right].	
\label{eq:R_def}
\end{equation}
For $n = 0$, the equation \eqref{eq:xi0} implies
\begin{equation*}
|\widehat{\xi}^0_l| \lesssim \varepsilon \int^{\tau}_0 \int^s_0 |\partial_{s_1} \widehat{F}^0_l(s_1)| d s_1 d s.
\end{equation*}
By the Parseval's identity and the assumptions (A) and (B), we find
\begin{align}
\|\textbf{e}^1(x)\|^2_{L^2} & = \|\xi^0(x)\|^2_{L^2} = (b-a) \sum_{l \in \mathcal{T}_N} |\widehat{\xi}^0_l|^2  \nonumber \\ 
& \lesssim \varepsilon^2 (b - a) \tau^2 \int^{\tau}_0 \int^s_0 \sum_{l \in \mathcal{T}_N} |\partial_{s_1}  \widehat{F}^0_l(s_1)|^2 d s_1 ds   \nonumber \\
& \lesssim \varepsilon^2 \tau^2  \int^{\tau}_0 \int^s_0 \|\partial_{s_1}  (G\Phi(s_1))\|^2_{L^2} d s_1 ds)\nonumber  \\
& \lesssim \varepsilon^2 \tau^4.	
\label{eq:e1_bound}
\end{align}
From \eqref{eq:xin}, we have for $1 \le n \le \frac{T_0/\varepsilon}{\tau}-1$,
\begin{equation*}
|\widehat{\xi}^n_l| \leq \varepsilon \int^{\tau}_0 \left(\int^s_0\int^{s_1}_0 |\partial_{s_2s_2} \widehat{F}^n_l(s_2)| d s_2 d s_1 + s\int^1_0 \int^{\tau}_{\theta\tau} |\partial_{\theta_1\theta_1} \widehat{F}^{n-1}_l(\theta_1)| d \theta_1 d\theta\right) d s.
\end{equation*}
Similar to the estimate \eqref{eq:e1_bound}, under the assumptions (A) and (B), we obtain
\begin{align}
\|\xi^n(x)\|^2_{L^2} = &\  (b-a) \sum_{l \in \mathcal{T}_N} |\widehat{\xi}^n_l|^2\nn \\
\lesssim & \  \varepsilon^2\tau^3 \int^{\tau}_0 \int^s_0\int^{s_1}_0 \sum_{l \in \mathcal{T}_N} |\partial_{s_2s_2} \widehat{F}^n_l(s_2)|^2 d s_2 d s_1 ds  \nn \\
& \ + \varepsilon^2\tau^3 \int^{\tau}_0 \int^1_0 \int^{\tau}_{\theta\tau} s \sum_{l \in \mathcal{T}_N}  |\partial_{\theta_1\theta_1} \widehat{F}^{n-1}_l(\theta_1)|^2 d \theta_1 d\theta d s  \nn  \\
\lesssim & \ \varepsilon^2\tau^6 \|\partial_{tt}F(t, x)\|^2_{L^{\infty}([0, T_0/\varepsilon];(L^2)^2)} \nn \\
 \lesssim & \  \varepsilon^2\tau^6, \quad 1 \le n \le \frac{T_0/\varepsilon}{\tau}-1.
\label{eq:xi_error}
\end{align}
Using the properties of the matrices $Q^{(1)}_l(\tau)$ and $Q^{(2)}_l(\tau)$, it is easy to check that
\begin{equation}
\|Q^{(1)}_l(\tau)\|_2 \leq \tau, \quad \|Q^{(2)}_l(\tau)\|_2\leq \frac{\tau^2}{2},\quad l \in \mathcal{T}_N.	
\label{eq:Q}
\end{equation}
Combining \eqref{eq:R_def} and \eqref{eq:Q}, we get for $1 \le n \le \frac{T_0/\varepsilon}{\tau}-1$,
\begin{align}
\|R^n\|^2_{L^2}  = & \ (b-a) \sum_{l \in \mathcal{T}_N}	|\widehat{R}^n_l|^2 \nonumber  \\
 \lesssim & \ \varepsilon^2 (b-a)\tau^2  \sum_{l \in \mathcal{T}_N}	\left[\left|\widehat{\left(G(t_n)\Phi(t_n)\right)}_l - \widehat{\left(G(t_n)(\Phi^n_N)\right)}_l\right|^2  + \left|\widehat{(G(t_{n-1})\Phi(t_{n-1}))}_l - \widehat{(G(t_{n-1})\Phi^{n-1}_N)}_l\right|^2\right] \nonumber \\
 \lesssim & \ \varepsilon^2 \tau^2 \left[\|G(t_n)\Phi(t_n, x) - G(t_n)\Phi^n_N(x)\|^2_{L^2} + \|G(t_{n-1})\Phi(t_{n-1}, x) - G(t_{n-1})\Phi^{n-1}_N(x)\|^2_{L^2}\right] \nonumber \\
 \lesssim & \ \varepsilon^2 \tau^2 h^{2m_0} + \varepsilon^2 \tau^2\|\textbf{e}^n(x)\|^2_{L^2} + \varepsilon^2 \tau^2\|\textbf{e}^{n-1}(x)\|^2_{L^2}.
\label{eq:R_error}
\end{align}
Multiplying both sides of \eqref{eq:errorfun} from left by $(\widehat{\textbf{e}}_l^{n+1} + e^{-i\tau\Gamma_l}\widehat{\textbf{e}}_l^{n})^{\ast}$, taking the real parts and using the Cauchy inequality, we obtain
\begin{equation}
|\widehat{\textbf{e}}^{n+1}_l|^2	- |\widehat{\textbf{e}}^{n}_l|^2 \leq \varepsilon \tau \left(|\widehat{\textbf{e}}^{n+1}_l|^2 + |\widehat{\textbf{e}}^{n}_l|^2\right) + \frac{|\widehat{R}^n_l|^2}{\varepsilon \tau}+ \frac{|\widehat{\xi}^n_l|^2}{\varepsilon \tau}.
\end{equation}
Summing up above inequalities for $l \in \mathcal{T}_N$ and then multiplying it by $(b - a)$,  by using the Parseval's identity, we get for $ n \geq 1$,
\begin{align}
\|\textbf{e}^	{n+1}(x)\|^2_{L^2} - \|\textbf{e}^	{n}(x)\|^2_{L^2}  \lesssim  \varepsilon\tau\left(\|\textbf{e}^	{n+1}(x)\|^2_{L^2} + \|\textbf{e}^{n}(x)\|^2_{L^2}\right)  + \frac{1}{\varepsilon \tau}\left(\|R^{n}(x)\|^2_{L^2} + \|\xi^n(x)\|^2_{L^2}\right). 
\label{eq:error1}
\end{align}
Summing up \eqref{eq:error1} for $ n = 1, \cdots, m$, combining \eqref{eq:xi_error} and \eqref{eq:R_error}, we obtain for $1 \leq m\leq \frac{T_0/\varepsilon}{\tau} - 1$,
\begin{equation}
\|\textbf{e}^{m+1}(x)\|^2_{L^2} - \|\textbf{e}^1(x)\|^2_{L^2} \lesssim \varepsilon \tau \sum^{m+1}_{k = 1} \|{\textbf e}^{k}(x)\|^2_{L^2}  + m\varepsilon\tau^5 + m\varepsilon\tau h^{2m_0}.
\end{equation}
Noticing $\|\textbf{e}^1(x)\|^2_{L^2} \lesssim \varepsilon^2 \tau^4$ and using the discrete Gronwall's inequality, there exists $0 < \tau_0 \leq \frac{1}{2}$ sufficiently small and independent of $\varepsilon$ such that for $0 < \varepsilon \leq 1$, when $0 < \tau \leq \tau_0$, we get
\begin{equation}
\|\textbf{e}^{m+1}(x)\|_{L^2} \leq (C_1 +C_2T_0)e^{C_3 T_0} \left( h^{m_0} + \tau^2 \right),\quad 1 \leq m\leq \frac{T_0/\varepsilon}{\tau} - 1,
\label{eq:C12}
\end{equation}
with $C_1$, $C_2$ and $C_3$ constants independent of $h, \tau$ and $\varepsilon$. Recalling $\|\textbf{e}^0(x)\|_{L^2}=0$ and $\|\textbf{e}^1(x)\|_{L^2} \lesssim \varepsilon\tau^2$, we have for $0\leq n\leq \frac{T_0/\varepsilon}{\tau}$,
\begin{equation*}
\|\Phi(t_{n}, x) - \Phi^{n}_N(x)\|_{L^2} \leq C_0 h^{m_0} + \|\textbf{e}^{n}(x)\|_{L^2} \leq  \left(C_0 + (C_1 + C_2T_0)e^{C_3 T_0}\right) \left( h^{m_0} + \tau^2 \right),
\end{equation*}
 which proves the desired error bound \eqref{eq:EWI_error} by taking $C_{T_0} = C_0 + (C_1 + C_2T_0)e^{C_3 T_0}$. The proof for Theorem \ref{thm:EWI-FS} is thus completed.

\subsection{Proof for Theorem \ref{thm:sEWI-FS}}
With the same notation in the proof for Theorem \ref{thm:EWI-FS}, we also just need to estimate $\|\textbf{e}^n(x)\|_{L^2}$ for $0 \leq n \leq \frac{T_0/\varepsilon}{\tau}$. Define the local truncation error $\eta^n(x) = \sum_{l \in \mathcal{T}_N} \widehat{\eta}^n_l e^{i \mu_l(x-a)} \in Y_N$ of the sEWI-FS \eqref{eq:sPhi_hat}-\eqref{eq:sPhi_hat_co} as
\begin{equation}
\begin{split}
\widehat{\eta}^0_l = &\ \widehat{(\Phi(\tau))}_l -  e^{-i\tau\Gamma_l}\widehat{(\Phi(0))}_l +  \varepsilon \Gamma_l^{-1}\left[I_2-e^{-i\tau\Gamma_l}\right] \widehat{(G(0)\Phi(0))}_l, \\
\widehat{\eta}^n_l = &\ \widehat{(\Phi(t_{n+1}))}_l - \widehat{(\Phi(t_{n-1}))}_l  + 2i\sin(\tau\Gamma_l)\widehat{(\Phi(t_n))}_l + i \frac{2\varepsilon}{\delta_l}\sin(\tau\delta_l)\widehat{(G(t_n)\Phi(t_n))}_l,\ n \geq 1.
\end{split}
\label{eq:local_sEWI}
\end{equation}
Similar to the analysis of the local truncation error for the EWI-FS method, it is easy to obtain
\begin{equation}
\|\eta^0(x)\|_{L^2} \lesssim \varepsilon \tau^2, \quad 	\|\eta^n(x)\|_{L^2} \lesssim \varepsilon \tau^3, \quad n \geq 1.
\label{eq:eta_error}
\end{equation}

Now, we are going to derive the error equations. Subtracting \eqref{eq:sPhi_hat_co} from \eqref{eq:local_sEWI}, we obtain the error equations as
\begin{align}
& \widehat{\mathbf{e}}^{n+1}_l -  \widehat{\mathbf{e}}^{n-1}_l = -2i\sin(\tau\Gamma_l)\widehat{\mathbf{e}}^{n}_l + \widehat{W}^n_l + \widehat{\eta}^n_l, \quad 1 \leq n \leq \frac{T_0/\varepsilon}{\tau} -  1,  \label{eq:errorfun_s} \\
& \widehat{\mathbf{e}}^0_l = \textbf{0}, \quad \widehat{\mathbf{e}}^1_l = \widehat{\eta}^0_l, \quad l \in \mathcal{T}_N, 	
\end{align}
where $W^n(x) = \sum_{l \in \mathcal{T}_N} \widehat{W}^n_l e^{i \mu_l(x-a)} \in Y_N$ for $ n \geq 1$, with $\widehat{W}^n_l$ given by 
\begin{equation}
\widehat{W}^n_l = \  {-i}\frac{2\varepsilon}{\delta_l} \sin(\tau\delta_l)\left(\widehat{(G(t_n)\Phi(t_n))}_l - \widehat{(G(t_n)\Phi^n_N)}_l\right).	
\label{eq:W_def}
\end{equation}
Since $|\sin(\tau\delta_l)/\delta_l| \leq \tau$, from \eqref{eq:W_def} and the assumption (B), we arrive at
\begin{align}
\|W^n(x)\|^2_{L^2}  = & \ (b-a) \sum_{l \in \mathcal{T}_N}	|\widehat{W}^n_l|^2 \nonumber  \\
 \lesssim & \ \varepsilon^2 (b-a)\tau^2  \sum_{l \in \mathcal{T}_N} \left|\widehat{\left(G(t_n)\Phi(t_n)\right)}_l - \widehat{\left(G(t_n)(\Phi^n_N)\right)}_l\right|^2  \nonumber \\
 \lesssim & \ \varepsilon^2 \tau^2 \|G(t_n)\Phi(t_n, x) - G(t_n)\Phi^n_N(x)\|^2_{L^2} \nonumber \\
 \lesssim & \ \varepsilon^2 \tau^2 h^{2m_0} + \varepsilon^2 \tau^2\|\textbf{e}^n(x)\|^2_{L^2}.
\label{eq:W_error}
\end{align}
Multiplying both sides of \eqref{eq:errorfun_s} from left by $(\widehat{\textbf{e}}_l^{n})^{\ast}$ and taking the real parts, we have
\begin{equation}
\text{Re}\left((\widehat{\textbf{e}}_l^{n})^{\ast} \widehat{\textbf{e}}_l^{n+1}\right) - \text{Re}\left((\widehat{\textbf{e}}_l^{n})^{\ast} \widehat{\textbf{e}}_l^{n-1}\right) = \text{Re}\left((\widehat{\textbf{e}}_l^{n})^{\ast} (\widehat{W}_l^{n+1}+ \widehat{\eta}_l^{n+1})\right),
\end{equation}
which implies
\begin{equation}
\left|\widehat{\textbf{e}}_l^{n+1}\right|^2 + \left|\widehat{\textbf{e}}_l^{n}\right|^2 - \left|\widehat{\textbf{e}}_l^{n+1}-\widehat{\textbf{e}}_l^{n}\right|^2 = \left|\widehat{\textbf{e}}_l^{n}\right|^2 + \left|\widehat{\textbf{e}}_l^{n-1}\right|^2 - \left|\widehat{\textbf{e}}_l^{n}-\widehat{\textbf{e}}_l^{n-1}\right|^2	 + 2 \text{Re}\left((\widehat{\textbf{e}}_l^{n})^{\ast} (\widehat{W}_l^{n+1}+ \widehat{\eta}_l^{n+1})\right).
\label{eq:se1}
\end{equation}
Multiplying both sides of \eqref{eq:errorfun_s} from left by $(\widehat{\textbf{e}}_l^{n+1} - 2\widehat{\textbf{e}}_l^{n} + \widehat{\textbf{e}}_l^{n-1})^{\ast}$ and taking the real parts, we obtain
\begin{align}
\left|\widehat{\textbf{e}}_l^{n+1}-\widehat{\textbf{e}}_l^{n}\right|^2 - \left|\widehat{\textbf{e}}_l^{n}-\widehat{\textbf{e}}_l^{n-1}\right|^2  = & \  2\text{Im}\left((\widehat{\textbf{e}}_l^{n+1})^{\ast} \sin(\tau\Gamma_l)\widehat{\textbf{e}}_l^{n} \right) - 2\text{Im}\left((\widehat{\textbf{e}}_l^{n})^{\ast} \sin(\tau\Gamma_l)\widehat{\textbf{e}}_l^{n-1} \right) \nonumber \\
&\  +  \text{Re}\left((\widehat{\textbf{e}}_l^{n+1}- 2\widehat{\textbf{e}}_l^{n} +\widehat{\textbf{e}}_l^{n-1})^{\ast} (\widehat{W}_l^{n}+ \widehat{\eta}_l^{n})\right).
\label{eq:se2}
\end{align}
By using Cauchy inequality, summing \eqref{eq:se1} and \eqref{eq:se2}, 
we arrive at
\begin{align}	
& \ \left|\widehat{\textbf{e}}_l^{n+1}\right|^2 + \left|\widehat{\textbf{e}}_l^{n}\right|^2 - 2\text{Im}\left((\widehat{\textbf{e}}_l^{n+1})^{\ast} \sin(\tau\Gamma_l)\widehat{\textbf{e}}_l^{n} \right) \nonumber \\
\leq & \  \left|\widehat{\textbf{e}}_l^{n}\right|^2 + \left|\widehat{\textbf{e}}_l^{n-1}\right|^2 - 2\text{Im}\left((\widehat{\textbf{e}}_l^{n})^{\ast} \sin(\tau\Gamma_l)\widehat{\textbf{e}}_l^{n-1} \right) \nonumber \\
&  + \varepsilon\tau \left( \left|\widehat{\textbf{e}}_l^{n+1}\right|^2 + \left|\widehat{\textbf{e}}_l^{n-1}\right|^2 \right) + \frac{1}{\varepsilon\tau}\left(\left|\widehat{W}_l^{n}\right|^2 + \left|\widehat{\eta}_l^{n}\right|^2\right).
\label{eq:se3}
\end{align}
Denote 
\begin{equation}
\mathcal{E}^n	 = \left\|{\textbf{e}}^{n+1}(x)\right\|^2_{L^2} + \left\|{\textbf{e}}^{n}(x)\right\|^2_{L^2} - 2(b-a) \sum_{l \in \mathcal{T}_N} \text{Im}\left((\widehat{\textbf{e}}_l^{n+1})^{\ast} \sin(\tau\Gamma_l)\widehat{\textbf{e}}_l^{n} \right),
\end{equation}
and it follows from the stability condition \eqref{eq:sEWI_con} that $|\sin(\tau\Gamma_l)\widehat{\mathbf{e}}_l^n| \leq |\sin(\tau \delta_l)||\widehat{\mathbf{e}}_l^n| \leq \sin(\pi/3)|\widehat{\mathbf{e}}_l^n| = \frac{\sqrt{3}}{2}|\widehat{\mathbf{e}}_l^n|$, which yields
\begin{align}
\mathcal{E}^n	 & \geq  \left\|{\textbf{e}}^{n+1}(x)\right\|^2_{L^2} +   \left\|{\textbf{e}}^{n}(x)\right\|^2_{L^2} - \frac{\sqrt{3}}{2} (b-a)   \sum_{l \in \mathcal{T}_N} \left(\left|\widehat{\textbf{e}}_l^{n+1}\right|^2 +  \left|\widehat{\textbf{e}}_l^{n}\right|^2\right)\nonumber \\
& = \frac{2-\sqrt{3}}{2} \left( \left\|{\textbf{e}}^{n+1}(x)\right\|^2_{L^2} + \left\|{\textbf{e}}^{n}(x)\right\|^2_{L^2} \right).
\label{eq:Eb}
\end{align}
Multiplying \eqref{eq:se3} by $(b-a)$ and summing up for $l \in \mathcal{T}_N$, we get for $n \geq1$
\begin{equation}
\mathcal{E}^n - \mathcal{E}^{n-1} \leq \varepsilon\tau\left( \left\|{\textbf{e}}^{n+1}(x)\right\|^2_{L^2} + \left\|{\textbf{e}}^{n}(x)\right\|^2_{L^2} + \left\|{\textbf{e}}^{n-1}(x)\right\|^2_{L^2} \right) + \frac{1}{\varepsilon\tau} \left(\left\|{W}^{n}(x)\right\|^2_{L^2} + \left\|{\eta}^{n}(x)\right\|^2_{L^2}\right).
\label{eq:error2}
\end{equation}
Summing up \eqref{eq:error2} for $ n = 1, \cdots, m$, combining \eqref{eq:eta_error} and \eqref{eq:W_error}, we obtain for $0 \leq m\leq \frac{T_0/\varepsilon}{\tau} - 1$,
\begin{equation}
\mathcal{E}^{m} - \mathcal{E}^{0} \lesssim \varepsilon \tau \sum^{m+1}_{k = 0} \|{\textbf e}^{k}(x)\|^2_{L^2}  + m\varepsilon\tau^5 + m\varepsilon\tau h^{2m_0}.
\end{equation}
Since ${\textbf e}^0(x) = {\textbf 0}$ and $\mathcal{E}^{m}$ is bounded from below \eqref{eq:Eb}, we have for $1 \leq m \leq \frac{T_0/\varepsilon}{\tau}-1$,
\begin{equation}
\frac{2-\sqrt{3}}{2} \left( \left\|{\textbf{e}}^{m+1}(x)\right\|^2_{L^2} + \left\|{\textbf{e}}^{m}(x)\right\|^2_{L^2} \right)	-  \left\|{\textbf{e}}^1(x)\right\|^2_{L^2} \lesssim \varepsilon \tau \sum^{m+1}_{k = 0} \|{\textbf e}^{k}(x)\|^2_{L^2}  + m\varepsilon\tau^5 + m\varepsilon\tau h^{2m_0}.
\end{equation}
Noticing $\|\textbf{e}^1(x)\|^2_{L^2} \lesssim \varepsilon^2 \tau^4$ and using the discrete Gronwall's inequality, there exists $0 < \tau_0 \leq \frac{1}{2}$ sufficiently small and independent of $\varepsilon$ such that for $0 < \varepsilon \leq 1$, when $0 < \tau \leq \tau_0$, we get
\begin{equation}
\|\textbf{e}^{m+1}(x)\|_{L^2} \leq (C_1+C_2T_0)e^{C_3 T_0} \left( h^{m_0} + \tau^2 \right),\quad 1 \leq m\leq \frac{T_0/\varepsilon}{\tau} - 1,
\label{eq:C12}
\end{equation}
with $C_1$, $C_2$ and $C_3$ constants independent of $h, \tau$ and $\varepsilon$. Recalling $\|\textbf{e}^0(x)\|_{L^2}=0$ and $\|\textbf{e}^1(x)\|_{L^2} \lesssim \varepsilon\tau^2$, we have for $0\leq n\leq \frac{T_0/\varepsilon}{\tau}$,
\begin{equation*}
\|\Phi(t_{n}, x) - \Phi^{n}_N(x)\|_{L^2} \leq C_0 h^{m_0} + \|\textbf{e}^{n}(x)\|_{L^2} \leq  \left(C_0 + (C_1 + C_2T_0)e^{C_3 T_0}\right) \left( h^{m_0} + \tau^2 \right),
\end{equation*}
 which proves the desired error bound \eqref{eq:sEWI_error} by taking $C_{T_0} = C_0 + (C_1 + C_2T_0)e^{C_3 T_0}$. The proof for Theorem \ref{thm:sEWI-FS} is thus completed.

\section{Numerical results}
In this section, we compare the accuracy of the EWI-FP and sEWI-FP methods for the long-time dynamics of the Dirac equation in 1D. Then we simulate the dynamics of the Dirac equation \eqref{eq:Dirac_21} in 2D with a honeycomb lattice potential.

\subsection{Convergence test and comparison}
In the numerical experiment, the problem is solved numerically on the domain $\Omega = (0, 2\pi)$ with periodic boundary conditions. We choose the electromagnetic potentials  as
\begin{equation*}
V(t, x) = \frac{2}{2+\cos(x)},\quad A_1(t, x) = \frac{1}{2+\cos(x)}, \quad x \in \Omega, \quad t \geq 0,
\end{equation*}
and the initial data as
\begin{equation*}
\phi_1(0, x) = \frac{1}{2+\cos(x)},\quad \phi_2(0, x) = \frac{1}{1+\sin(x)^2},\quad x \in \Omega.
\end{equation*}
Since the exact solution is unknown, we use the TSFP method with a fine mesh size $h_e = \pi/64$  and a very small time step $\tau_e = 10^{-5}$ to get the `reference' solution numerically. Denote $\Phi^n_{h, \tau}$ as the numerical solution obtained by a numerical method with the mesh size $h$ and time step $\tau$. In order to quantify the numerical errors, we introduce the following discrete $l^2$ errors of the wave function $\Phi$
\begin{equation*}
e_{h, \tau}(t_n) = \|\Phi^n_{h, \tau} - \Phi(t_n, \cdot)\|_{l^2} = \sqrt{h\sum^{N-1}_{j=0} |\Phi^n_j - \Phi(t_n, x_j)|^2}.
\end{equation*}

\begin{table}[h!]
\renewcommand{\arraystretch}{1.3}
\def\temptablewidth{1\textwidth}
\setlength{\tabcolsep}{6pt}
\caption{Spatial errors of the EWI-FP method for the Dirac equation \eqref{eq:Dirac_1D} at $t=2/\varepsilon$.}
\label{tab:EWIFP_h}
{\rule{\temptablewidth}{1pt}}
\centering
\begin{tabular*}{\temptablewidth}{@{\extracolsep{\fill}}cccccc}
$e_{h,\tau_e}(t=2/\varepsilon)$ &$h_0 = \frac{\pi}{4}$ & $h_0/2 $ &$h_0/2^2 $ & $h_0/2^3$  \\
\hline
$\varepsilon_0 = 1$ & 1.27E-1 & 2.76E-3 & 3.34E-6 & 1.41E-9 \\
\hline
$\varepsilon_0 / 2$ & 1.25E-1 & 2.55E-3 & 2.90E-6 & 6.34E-10 \\
\hline
$\varepsilon_0 / 2^2$ & 1.09E-1 & 3.07E-3 & 1.29E-6 & 4.01E-10 \\
\hline
$\varepsilon_0 / 2^3$ & 2.77E-2 & 8.52E-4 & 2.46E-6 & 4.00E-10 \\
\hline
$\varepsilon_0 / 2^4$ & 3.06E-2 & 1.28E-3 & 2.19E-6 & 4.75E-10 \\
\hline
$\varepsilon_0 / 2^5$ & 3.08E-2 & 2.24E-3 & 2.63E-6 & 6.45E-10 \\
\end{tabular*}
{\rule{\temptablewidth}{1pt}}
\end{table}

\begin{figure}[ht!]
\begin{minipage}{0.5\textwidth}
\centerline{\includegraphics[width=8cm,height=7cm]{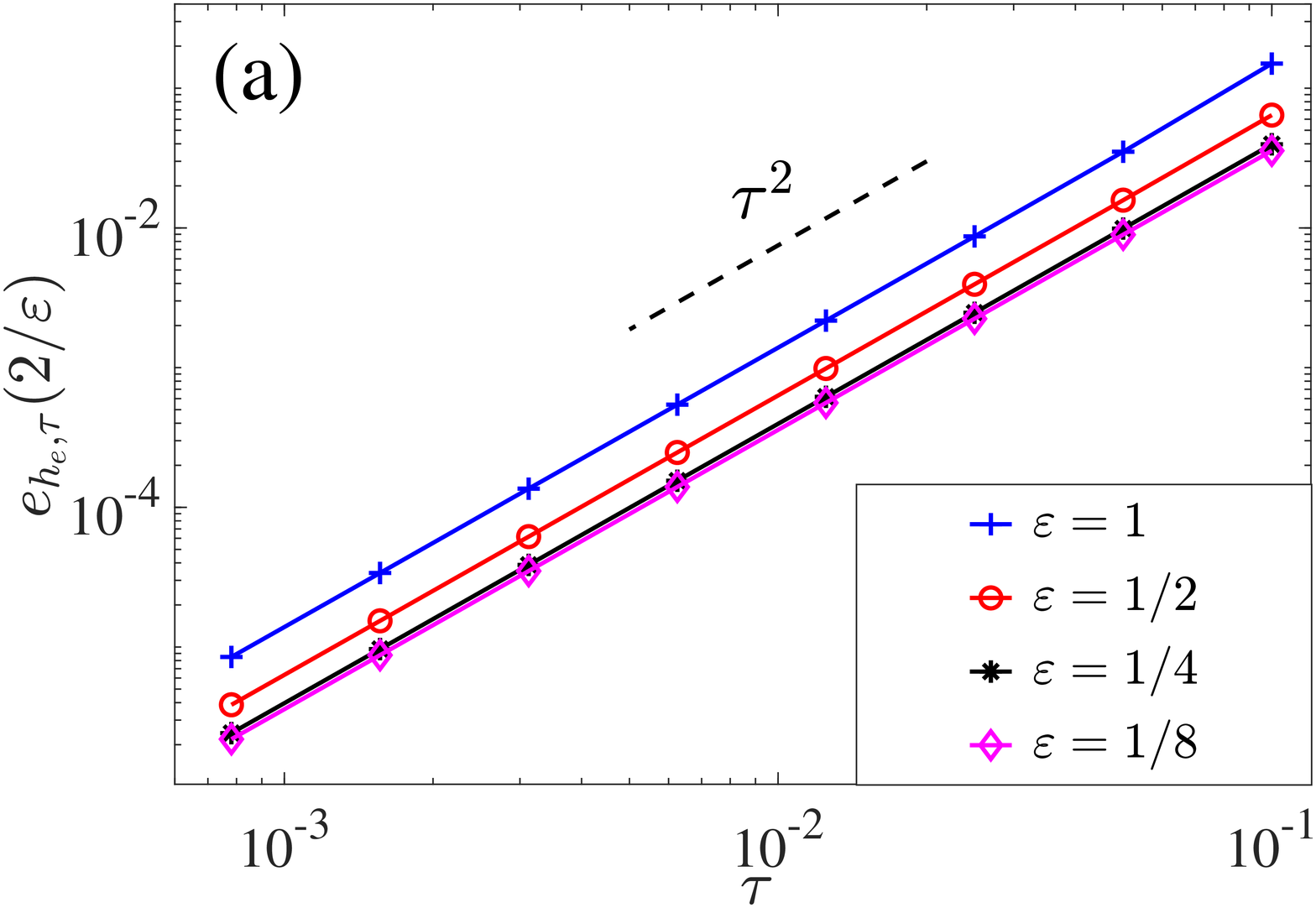}}
\end{minipage}
\begin{minipage}{0.5\textwidth}
\centerline{\includegraphics[width=8cm,height=7cm]{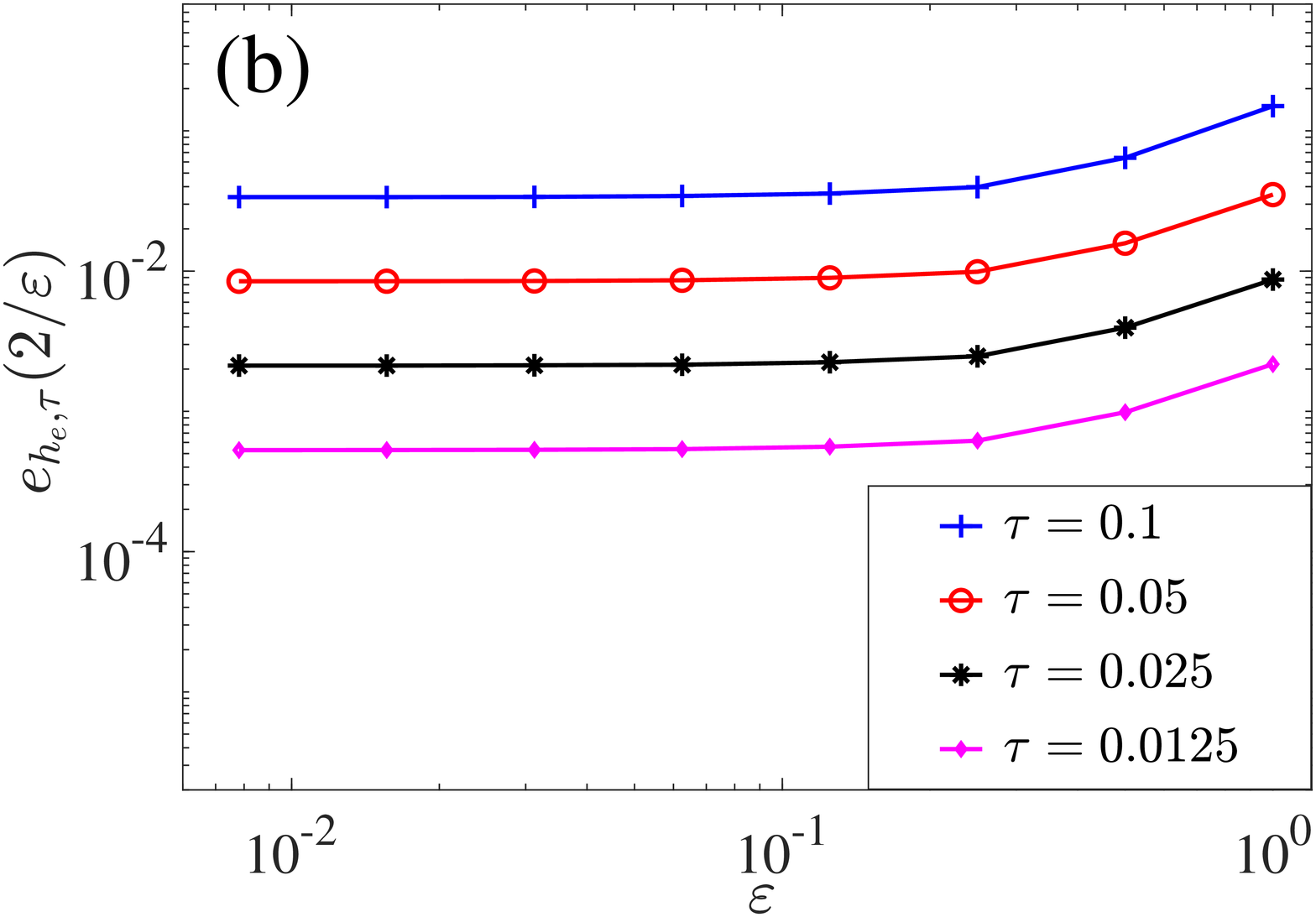}}
\end{minipage}
\caption{Temporal errors of the EWI-FP method for the Dirac equation \eqref{eq:Dirac_1D} at $t = 2/\varepsilon$.}
\label{fig:EWI_t}
\end{figure}

\begin{figure}[ht!]
\begin{minipage}{0.5\textwidth}
\centerline{\includegraphics[width=8cm,height=7cm]{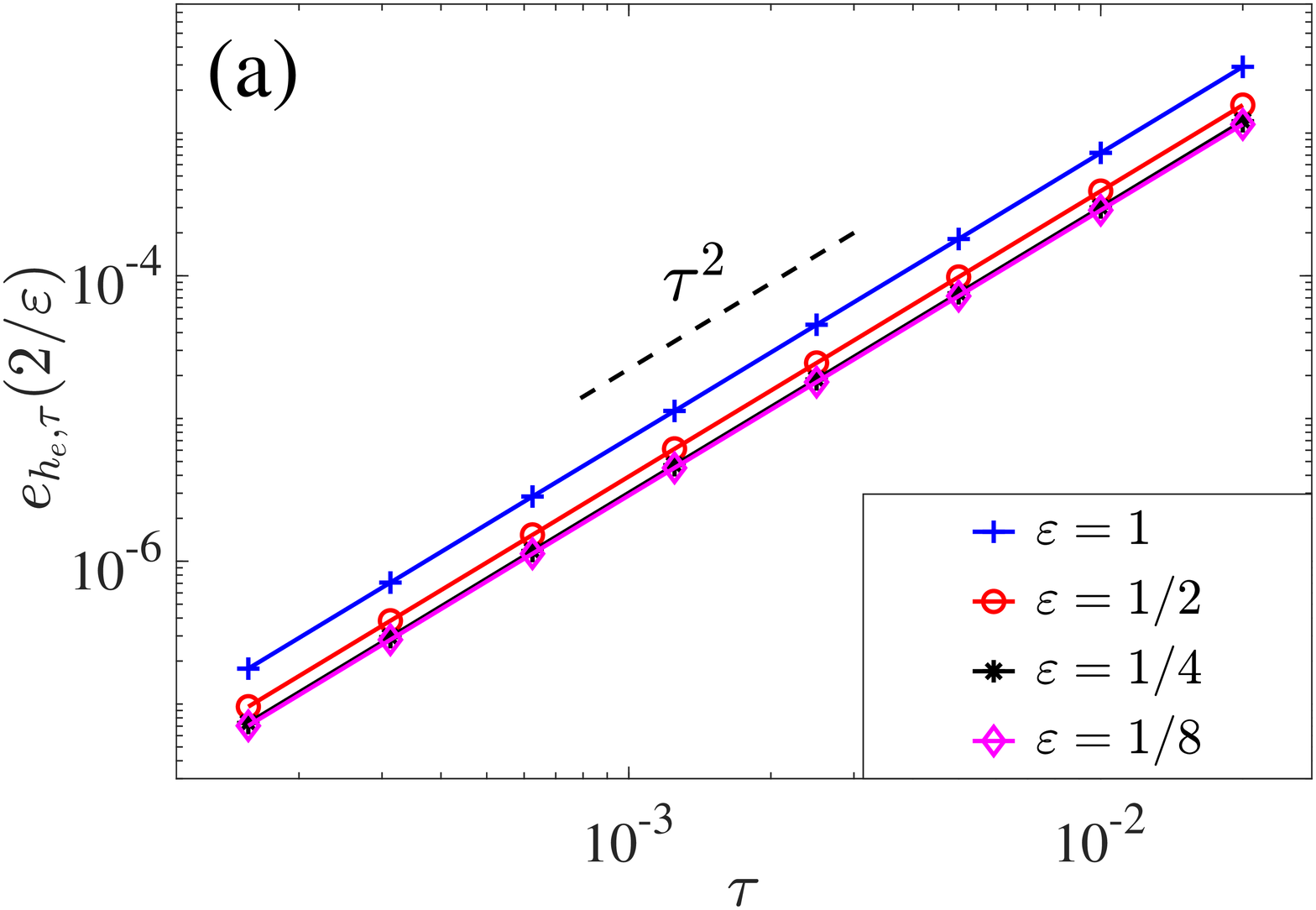}}
\end{minipage}
\begin{minipage}{0.5\textwidth}
\centerline{\includegraphics[width=8cm,height=7cm]{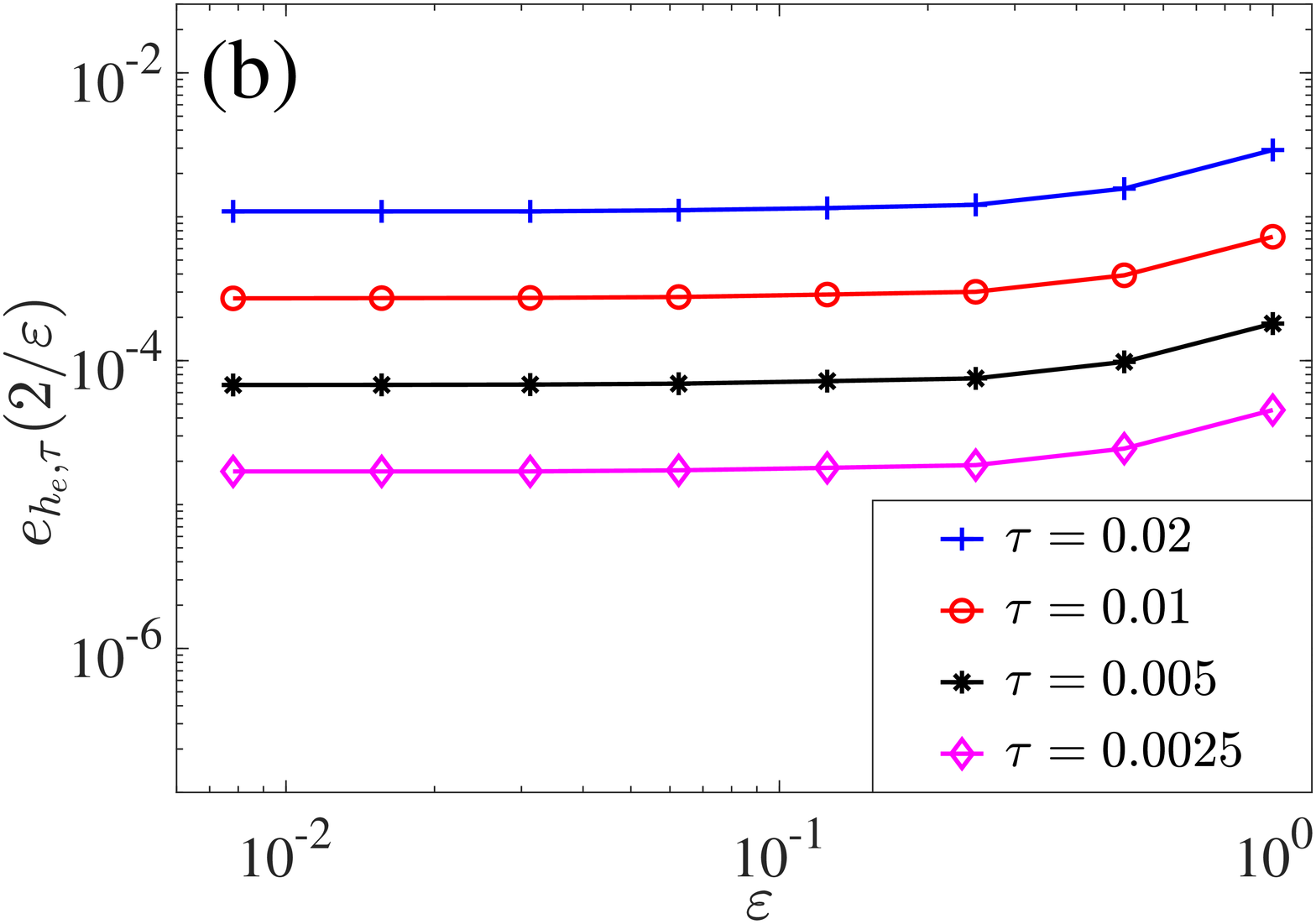}}
\end{minipage}
\caption{Temporal errors of the sEWI-FP method for the Dirac equation \eqref{eq:Dirac_1D} at $t = 2/\varepsilon$.}
\label{fig:sEWI_t}
\end{figure}

\begin{figure}[ht!]
\begin{minipage}{0.5\textwidth}
\centerline{\includegraphics[width=8cm,height=7cm]{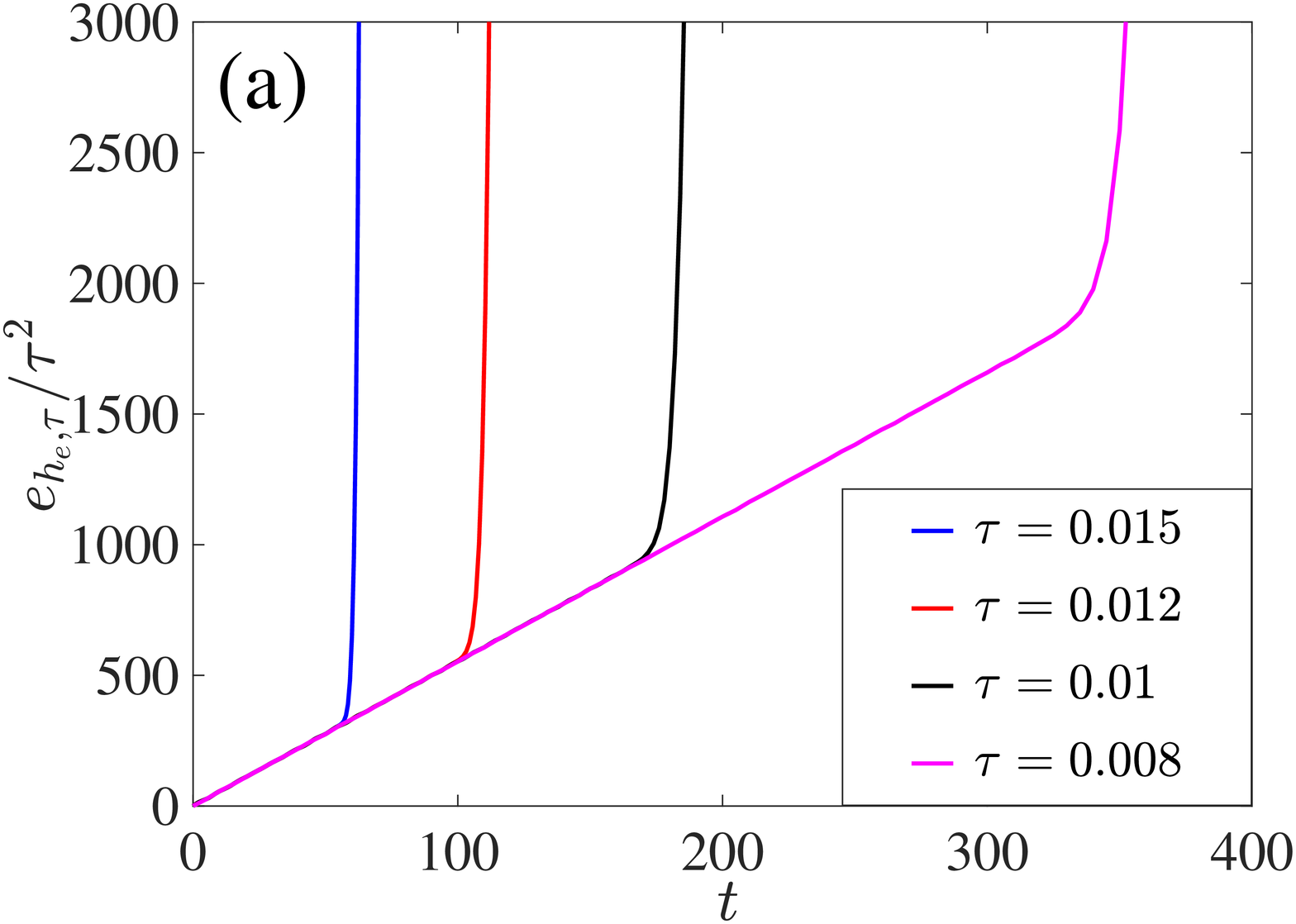}}
\end{minipage}
\begin{minipage}{0.5\textwidth}
\centerline{\includegraphics[width=8cm,height=7cm]{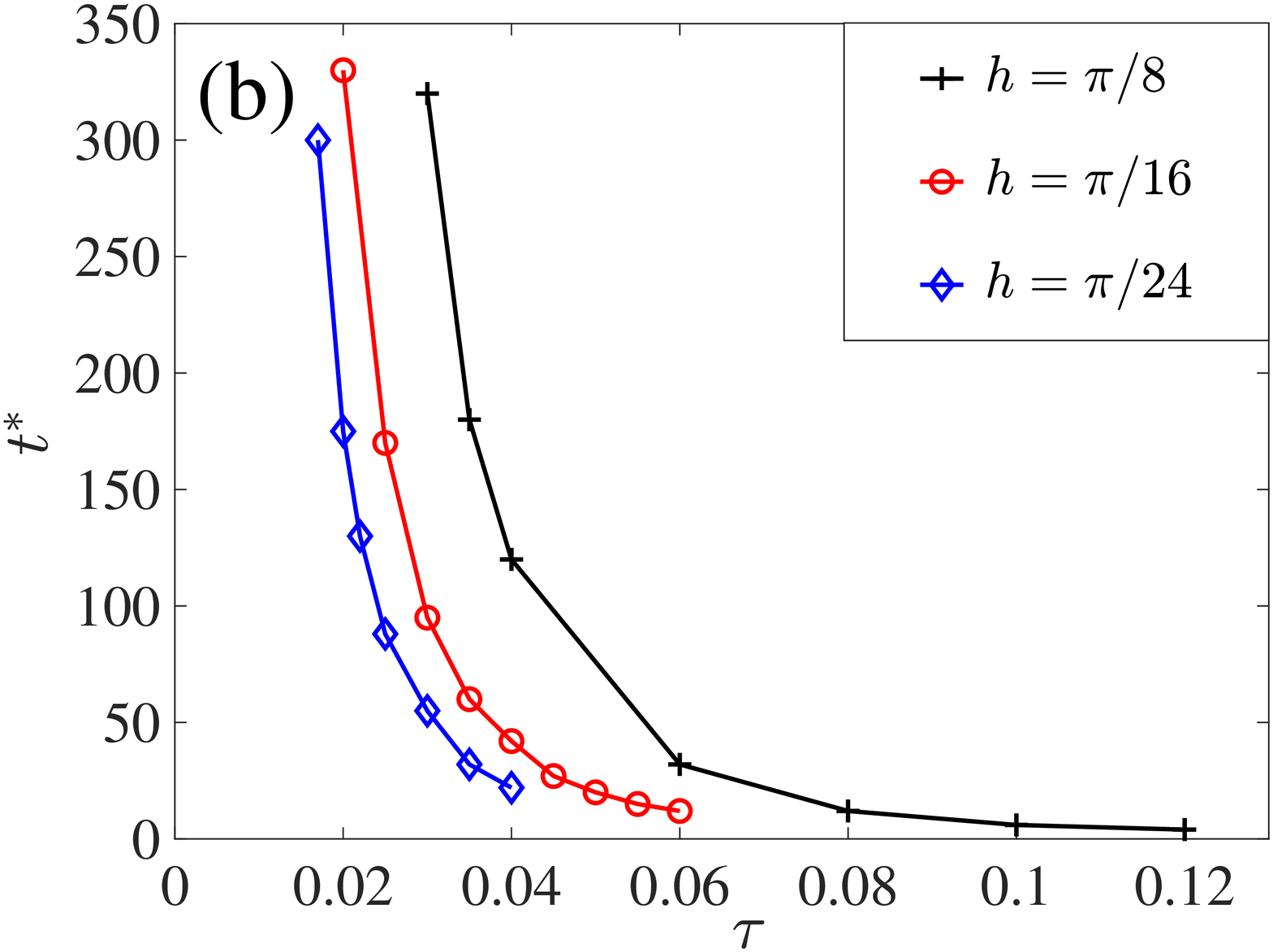}}
\end{minipage}
\caption{Long-time behaviors of the EWI-FP method for the Dirac equation \eqref{eq:Dirac_1D}  with $\varepsilon = 1$.}
\label{fig:EWI_temporal}
\end{figure}

\begin{figure}[ht!]
\centerline{\includegraphics[width=16cm,height=6cm]{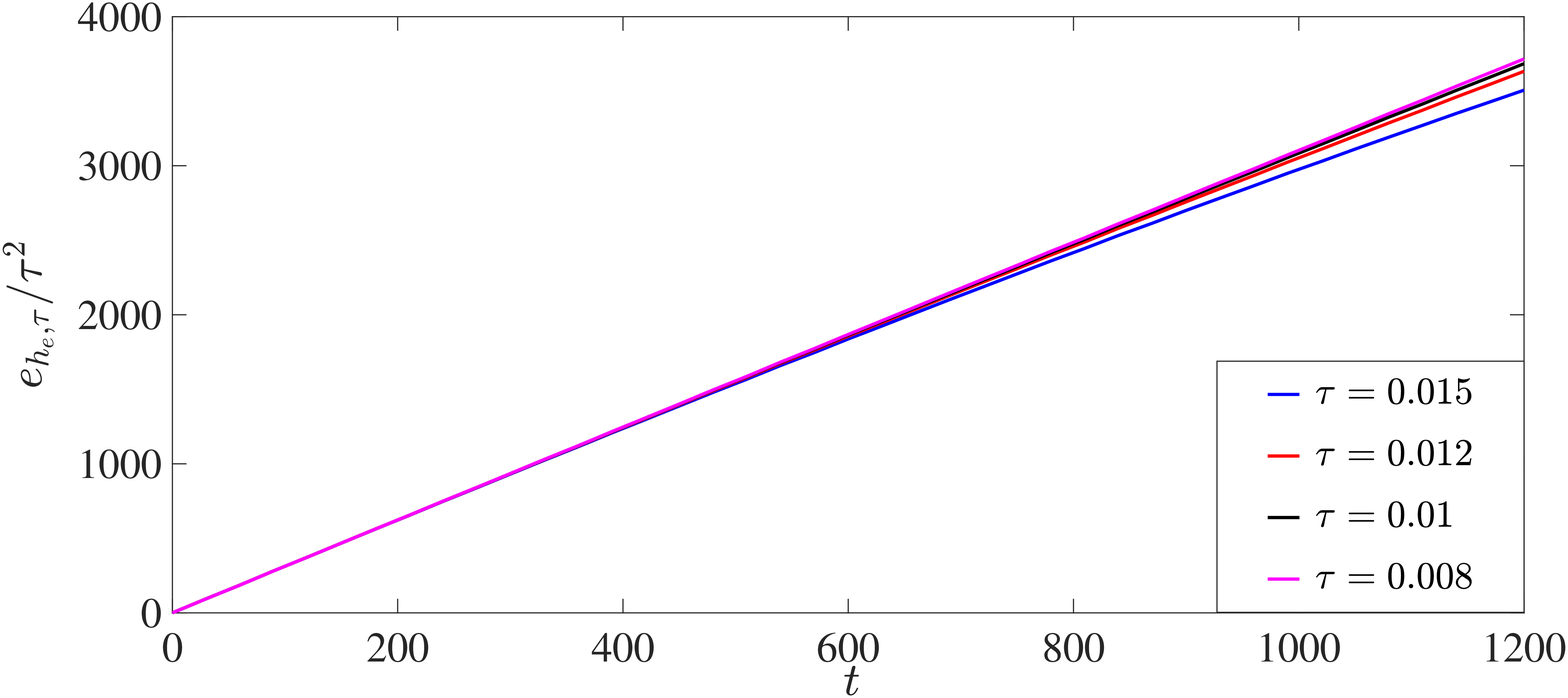}}
\caption{Long-time behaviors of the sEWI-FP method for the Dirac equation \eqref{eq:Dirac_1D} with $\varepsilon = 1$.}
\label{fig:sEWI_temporal}
\end{figure}

\begin{figure}[ht!]
\centerline{\includegraphics[width=16cm,height=6cm]{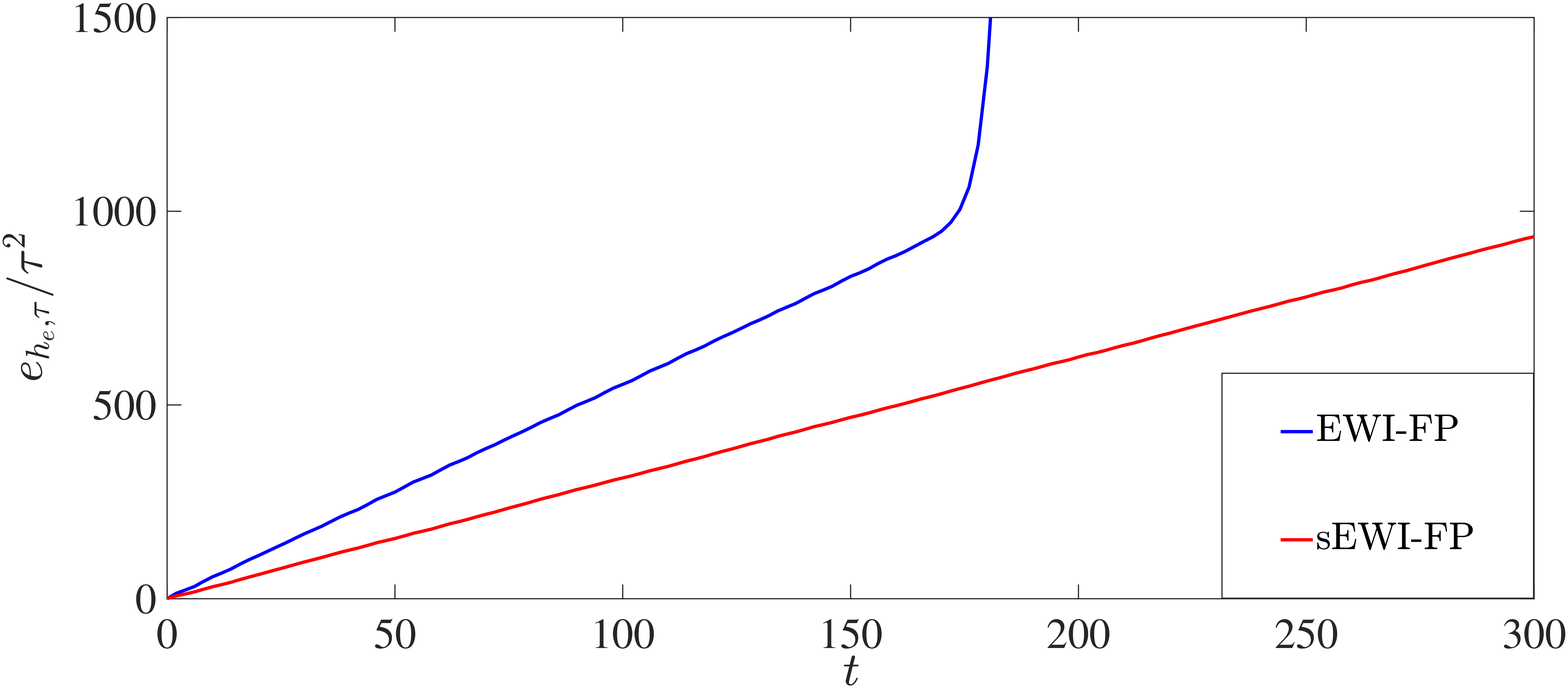}}
\caption{Comparison of the EWI-FP and sEWI-FP methods for the Dirac equation \eqref{eq:Dirac_1D} with $\varepsilon = 1$.}
\label{fig:com}
\end{figure}

Table \ref{tab:EWIFP_h} shows spatial errors $e_{h, \tau_e}(t=2/\varepsilon)$ of the EWI-FP method for different $h$ and $\varepsilon$ with the small time step $\tau_e = 10^{-5}$ such that the temporal discretization errors are negligible. The spatial errors of the sEWI-FP method are similar to the EWI-FP method and we omit here for brevity. Figure \ref{fig:EWI_t} and Figure \ref{fig:sEWI_t} depict temporal errors $e_{h_e, \tau}(t=2/\varepsilon)$ of the EWI-FP and sEWI-FP methods for different $\tau$ and $\varepsilon$, respectively. 

From Table \ref{tab:EWIFP_h}, Figures \ref{fig:EWI_t}-\ref{fig:sEWI_t}, we can draw the following conclusions for the long-time dynamics of the Dirac equation:

(i) The EWI-FP and sEWI-FP methods are both spectrally accurate in space (cf. each row in Table \ref{tab:EWIFP_h}) and the spatial errors are independent of $\varepsilon$ (cf. each column in Table \ref{tab:EWIFP_h}).

(ii) For any fixed $\varepsilon = \varepsilon_0 > 0$, the EWI-FP and sEWI-FP methods are both second-order in time (cf. Figure \ref{fig:EWI_t} (a) and Figure \ref{fig:sEWI_t} (a)). In the long-time regime, i.e., $0 < \varepsilon \ll 1$, the temporal errors of the EWI-FP and sEWI-FP methods are uniform for any $0 < \varepsilon \leq 1$ (cf. Figure \ref{fig:EWI_t} (b) and Figure \ref{fig:sEWI_t} (b)). 

Figure \ref{fig:EWI_temporal}(a) shows temporal errors of the EWI-FP method with different time step size $\tau$, which indicates that the temporal errors first grow linearly and then exponentially after the time $t^{\ast}$. Figure \ref{fig:EWI_temporal}(b) depicts the dependency of $t^{\ast}$ on the mesh size $h$ and time step size $\tau$. For the fixed mesh size $h$, the time $t^{\ast}$ is larger for the smaller time step size, which means that it could get longer-time simulations. For the fixed time step size $\tau$, the time $t^{\ast}$ is larger for the larger mesh size $h$. Figure \ref{fig:sEWI_temporal} shows the long-time behaviors of the sEWI-FP method for the Dirac equation \eqref{eq:Dirac_1D} with different time step size $\tau$  with $\varepsilon = 1$. Although we can only prove that the temporal error depends exponentially on the time $T_0$, the numerical simulations imply linear dependence, which is better than the analytical result. Figure \ref{fig:com} compares the long-time behaviors of the EWI-FP and sEWI-FP methods, and indicates that the symmetric scheme performs much better in the long-time simulations.

\subsection{Dynamics of the Dirac equation in 2D}
In this subsection, we study the dynamics of the Dirac equation \eqref{eq:Dirac_21} in 2D with a honeycomb lattice potential, i.e., we take $d = 2$, $A_1(t, {\bf x}) = A_2(t, {\bf x}) \equiv 0$ and 
\begin{equation}
V(t, {\bf x}) = \cos\left(\frac{4\pi}{\sqrt{3}}{\bf e}_1 \cdot {\bf x}\right)	+  \cos\left(\frac{4\pi}{\sqrt{3}}{\bf e}_2 \cdot {\bf x}\right)	+  \cos\left(\frac{4\pi}{\sqrt{3}}{\bf e}_3 \cdot {\bf x}\right),	
\end{equation}
with 
\begin{equation}
{\bf e}_1 = (-1, 0)^T, \quad {\bf e}_2 = (1/2, \sqrt{3}/2)^T, \quad  {\bf e}_3 = (1/2, -\sqrt{3}/2)^T.
\end{equation}
The initial data in \eqref{eq:initial} is chosen as 
\begin{equation}
\phi_1(0, {\bf x}) = e^{-\frac{x^2+y^2}{2}}, \quad \phi_2(0, {\bf x}) = e^{-\frac{x^2+y^2}{2}}, \quad {\bf x} = (x, y)^T \in \Omega = (-15, 15)^2.
\end{equation}
The problem is solved numerically by the EWI-FP/sEWI-FP method with the mesh size $h = 1/16$ and time step $\tau = 0.01$. 

\begin{figure}[ht!]
\begin{minipage}{0.5\textwidth}
\centerline{\includegraphics[width=8cm,height=6cm]{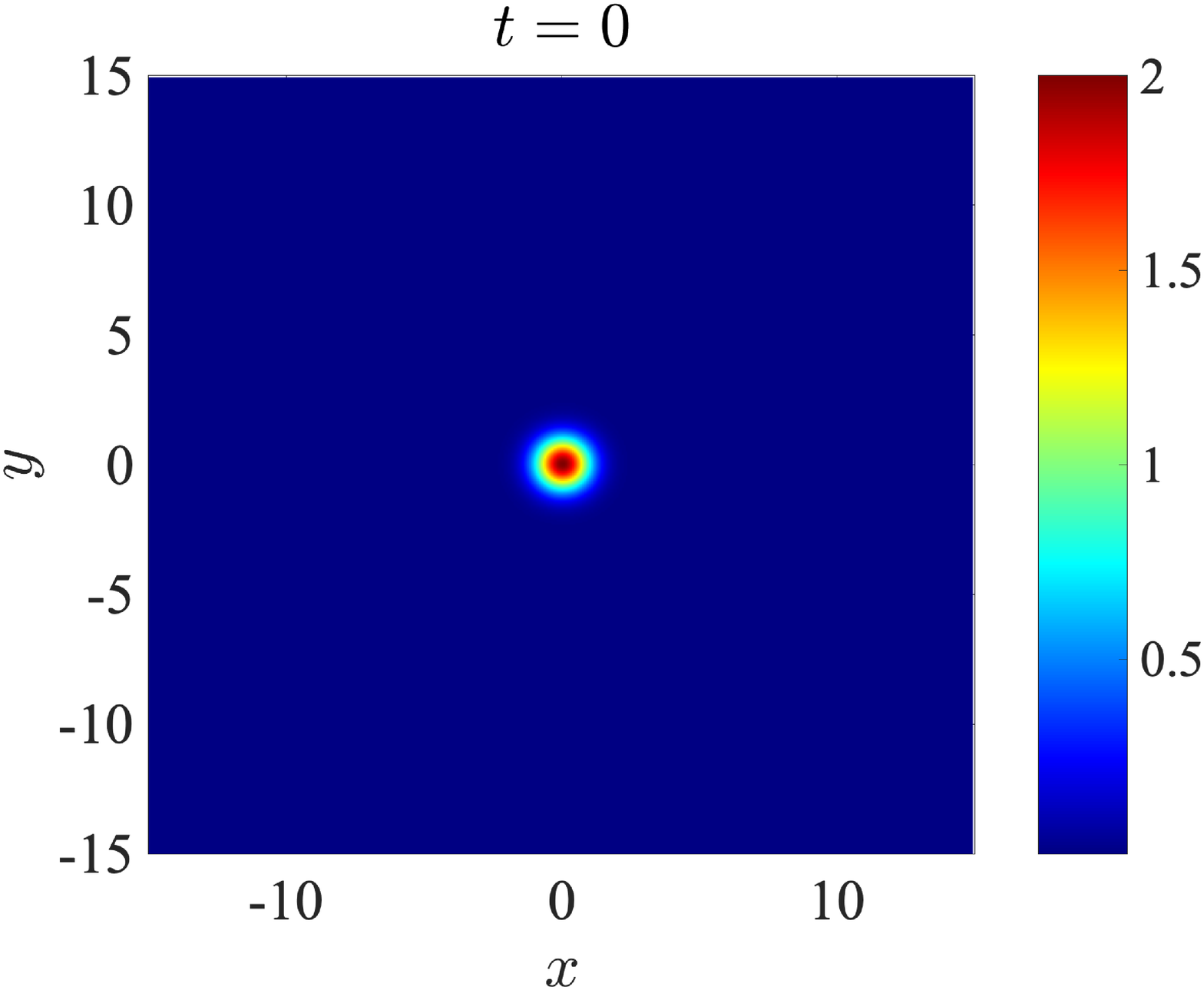}}
\end{minipage}
\begin{minipage}{0.5\textwidth}
\centerline{\includegraphics[width=8cm,height=6cm]{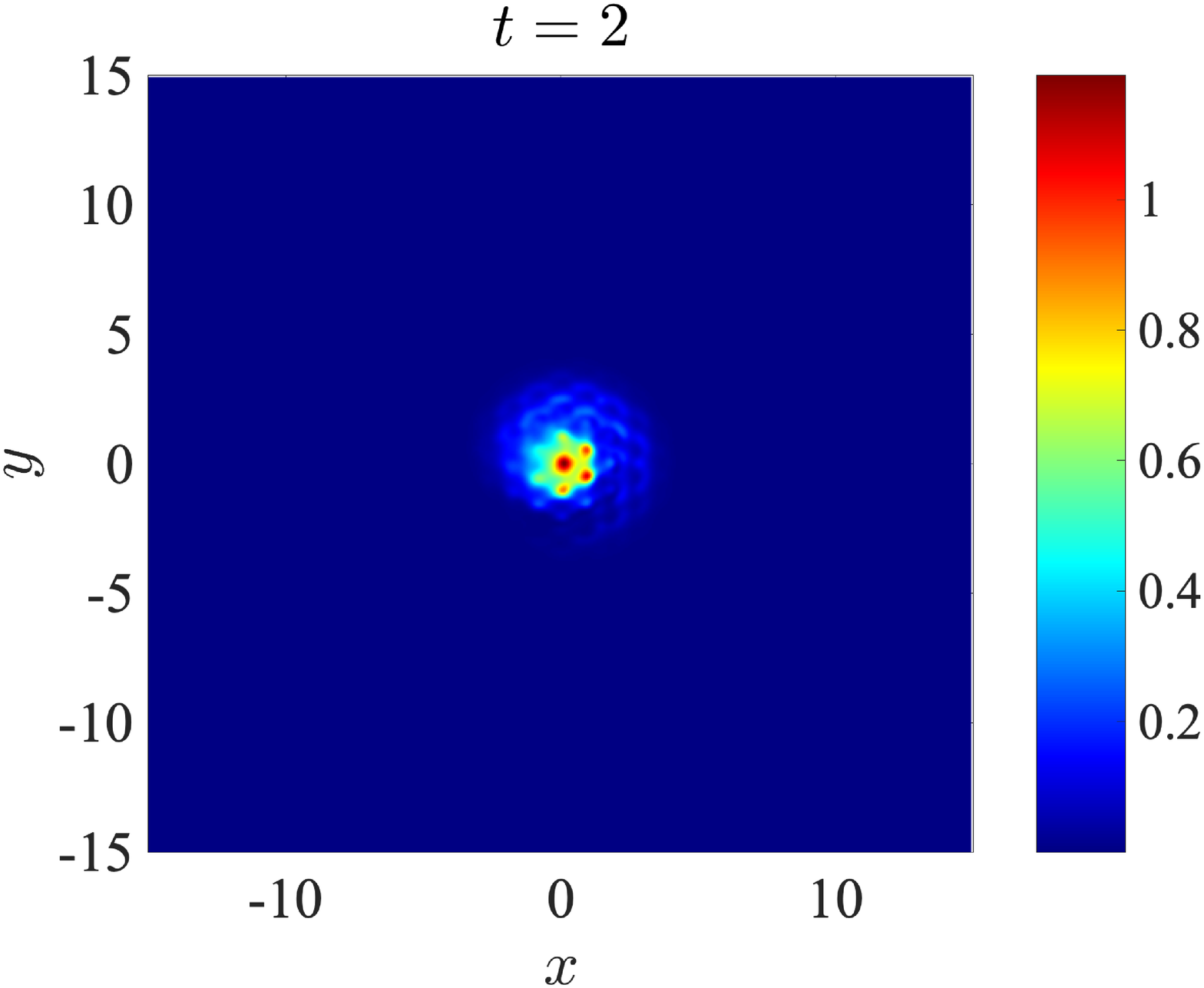}}
\end{minipage}
\begin{minipage}{0.5\textwidth}
\centerline{\includegraphics[width=8cm,height=6cm]{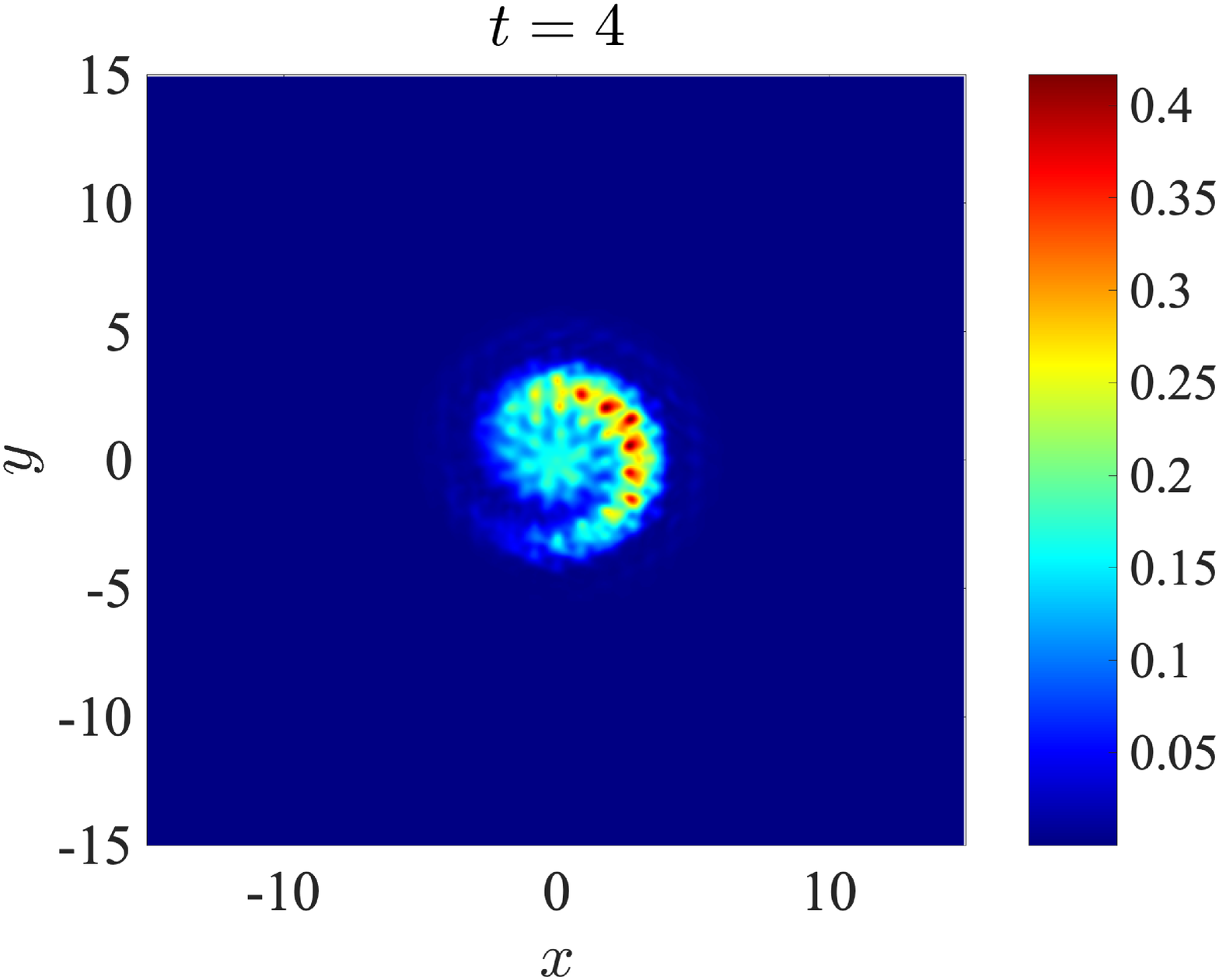}}
\end{minipage}
\begin{minipage}{0.5\textwidth}
\centerline{\includegraphics[width=8cm,height=6cm]{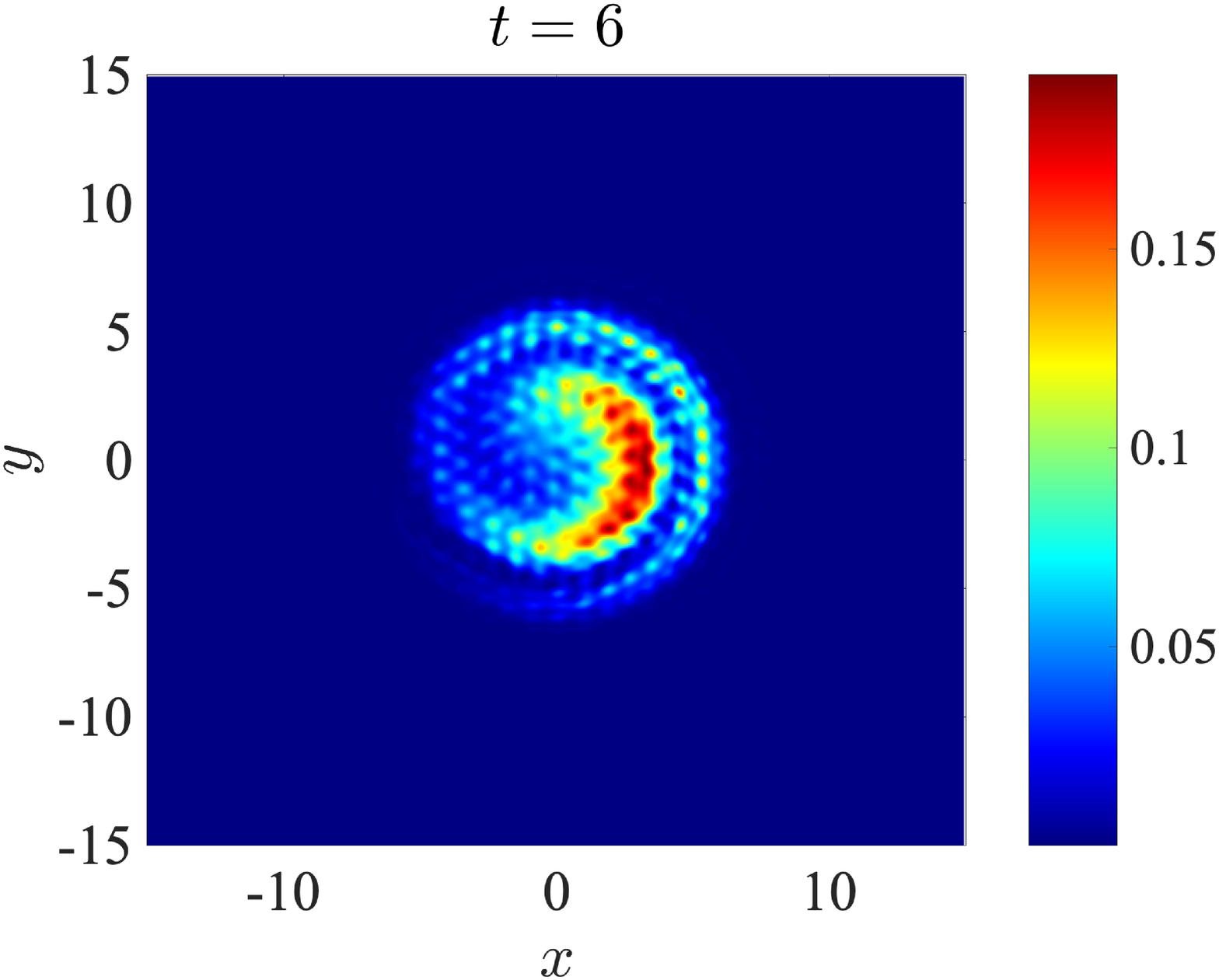}}
\end{minipage}
\begin{minipage}{0.5\textwidth}
\centerline{\includegraphics[width=8cm,height=6cm]{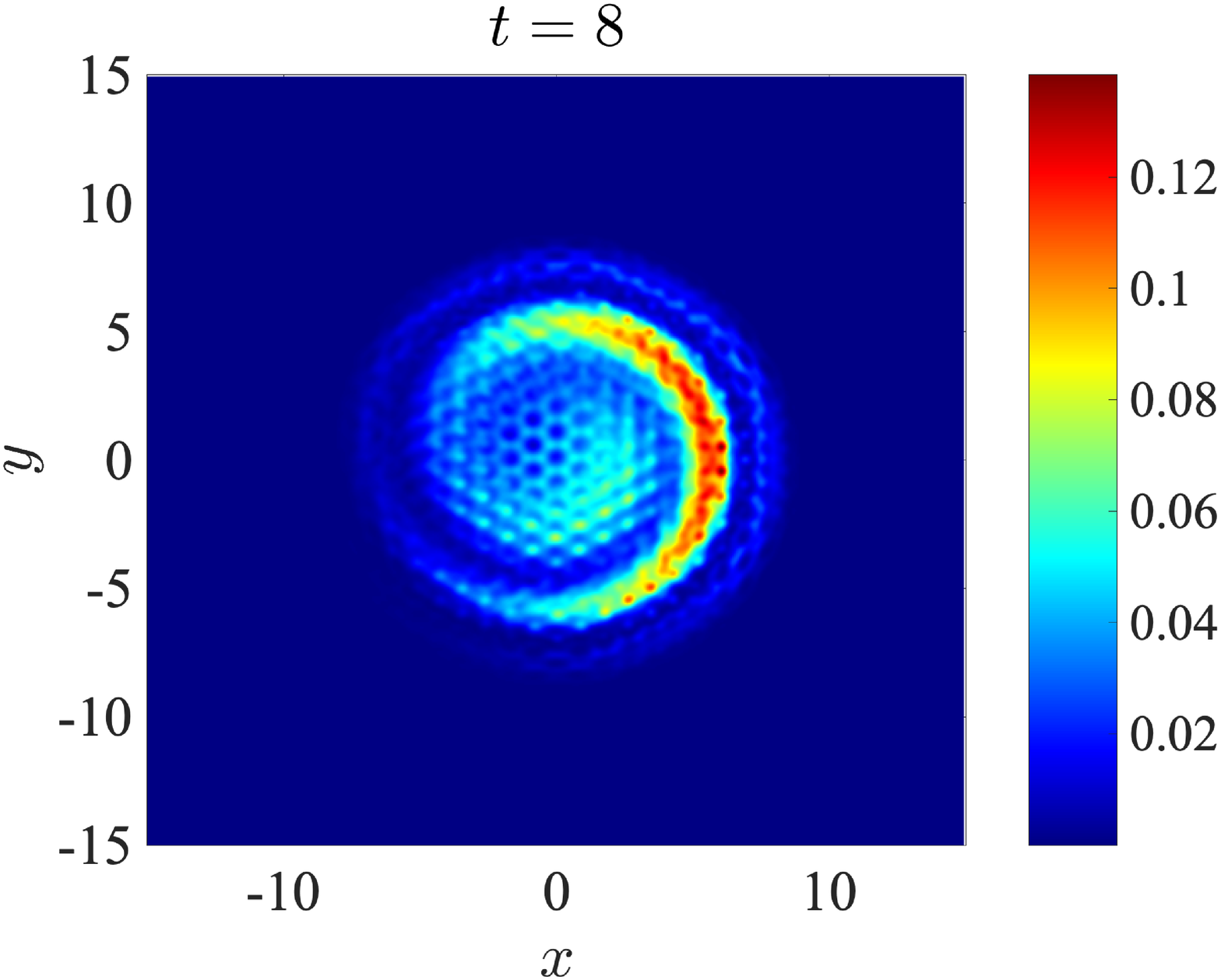}}
\end{minipage}
\begin{minipage}{0.5\textwidth}
\centerline{\includegraphics[width=8cm,height=6cm]{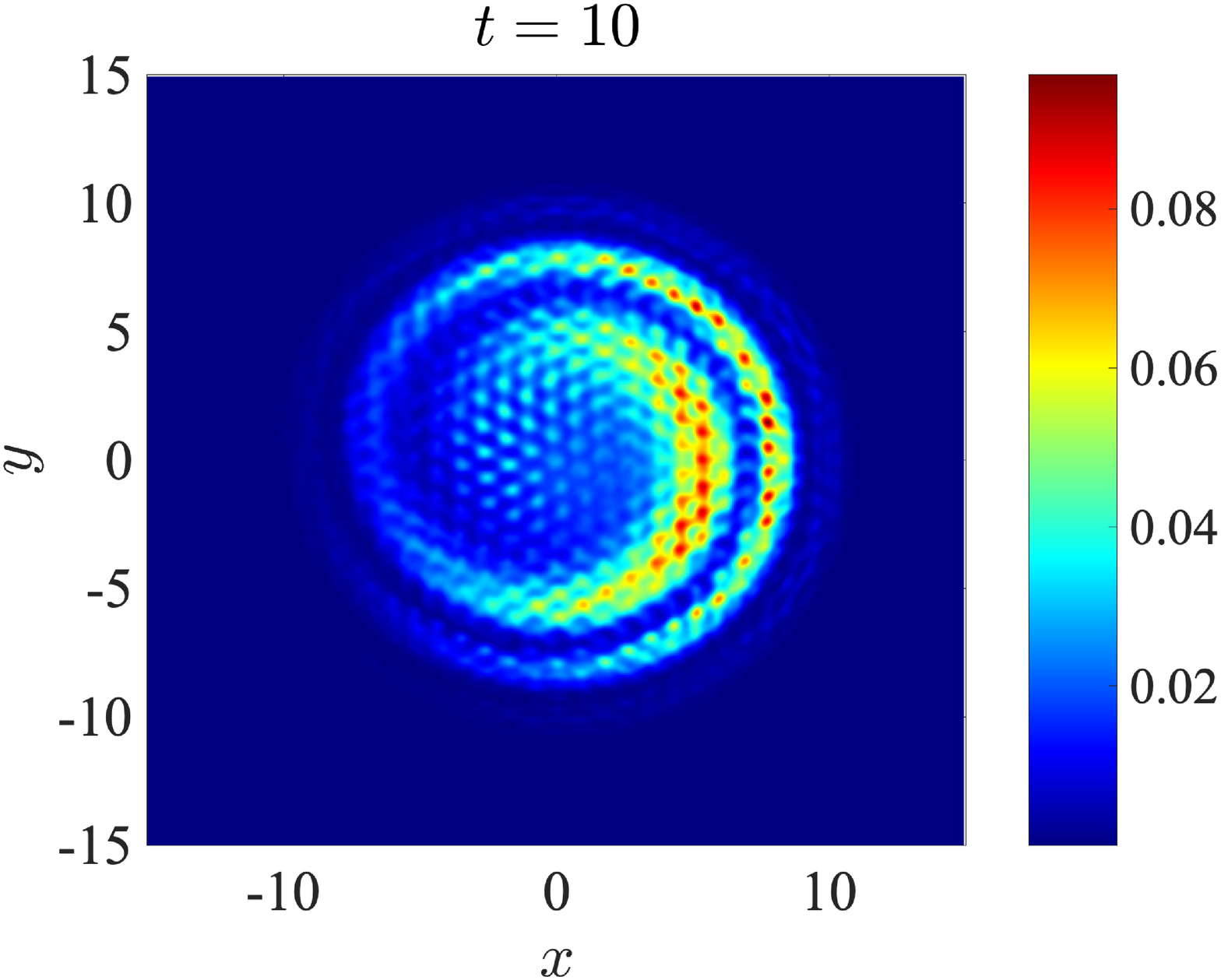}}
\end{minipage}
\caption{Dynamics of the density $\rho(t, {\bf x}) = |\phi_1(t, {\bf x})|^2 + |\phi_2(t, {\bf x})|^2 $ of the Dirac equation \eqref{eq:Dirac_21} in 2D with a honeycomb lattice potential when $\varepsilon = 1$ up to $t = 10$.}
\label{fig:eps1_rho}
\end{figure}

\begin{figure}[ht!]
\begin{minipage}{0.5\textwidth}
\centerline{\includegraphics[width=8cm,height=6cm]{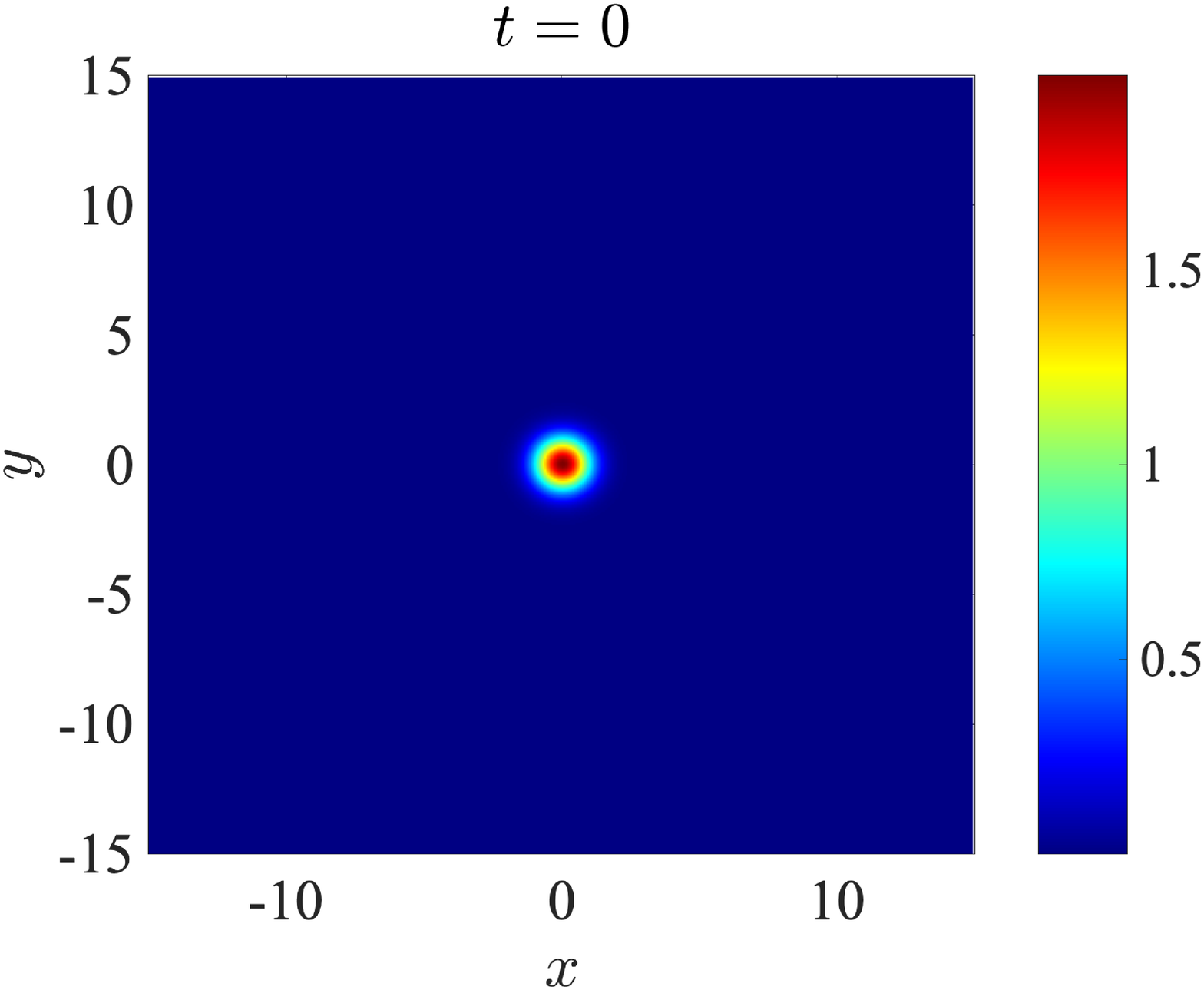}}
\end{minipage}
\begin{minipage}{0.5\textwidth}
\centerline{\includegraphics[width=8cm,height=6cm]{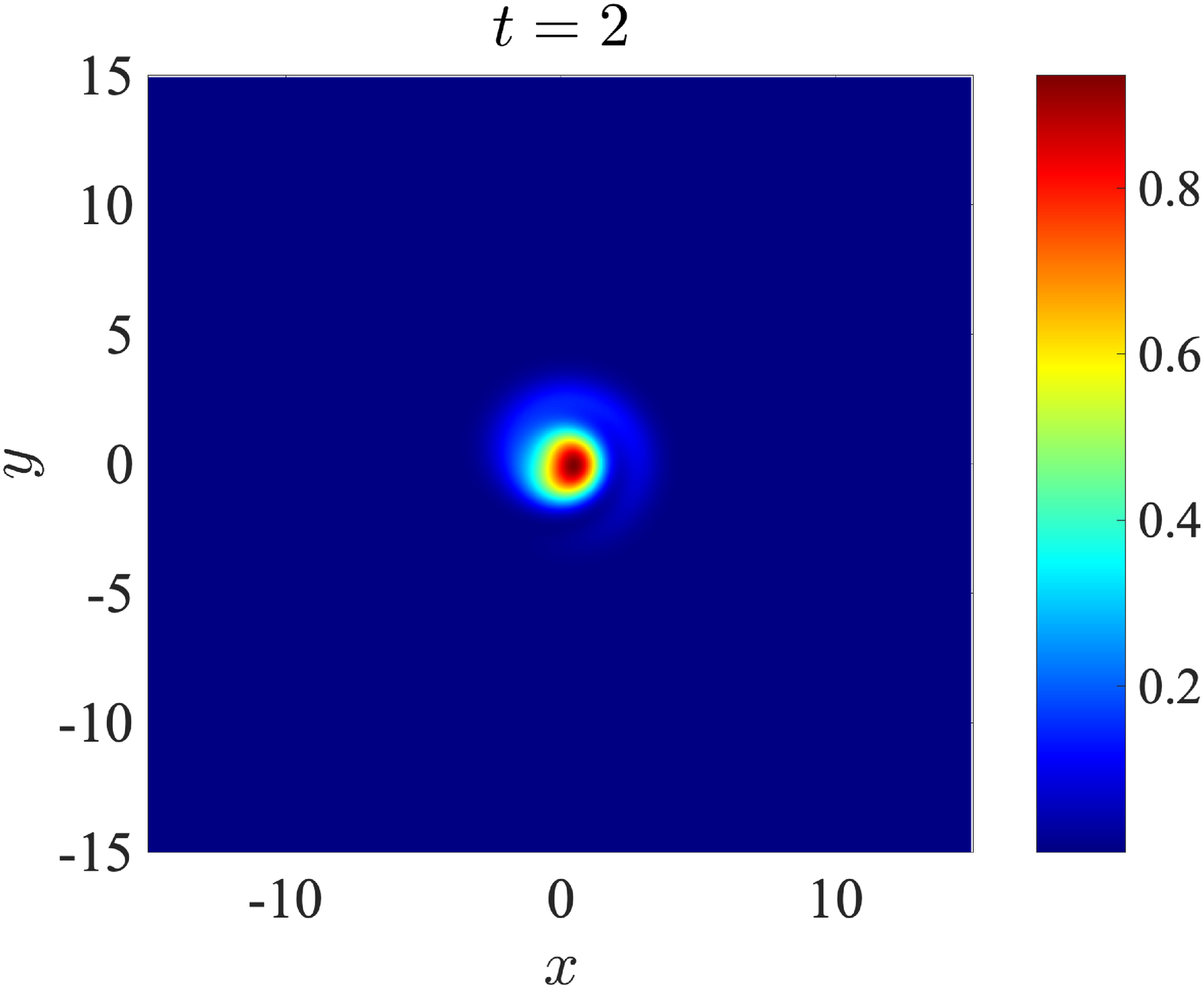}}
\end{minipage}
\begin{minipage}{0.5\textwidth}
\centerline{\includegraphics[width=8cm,height=6cm]{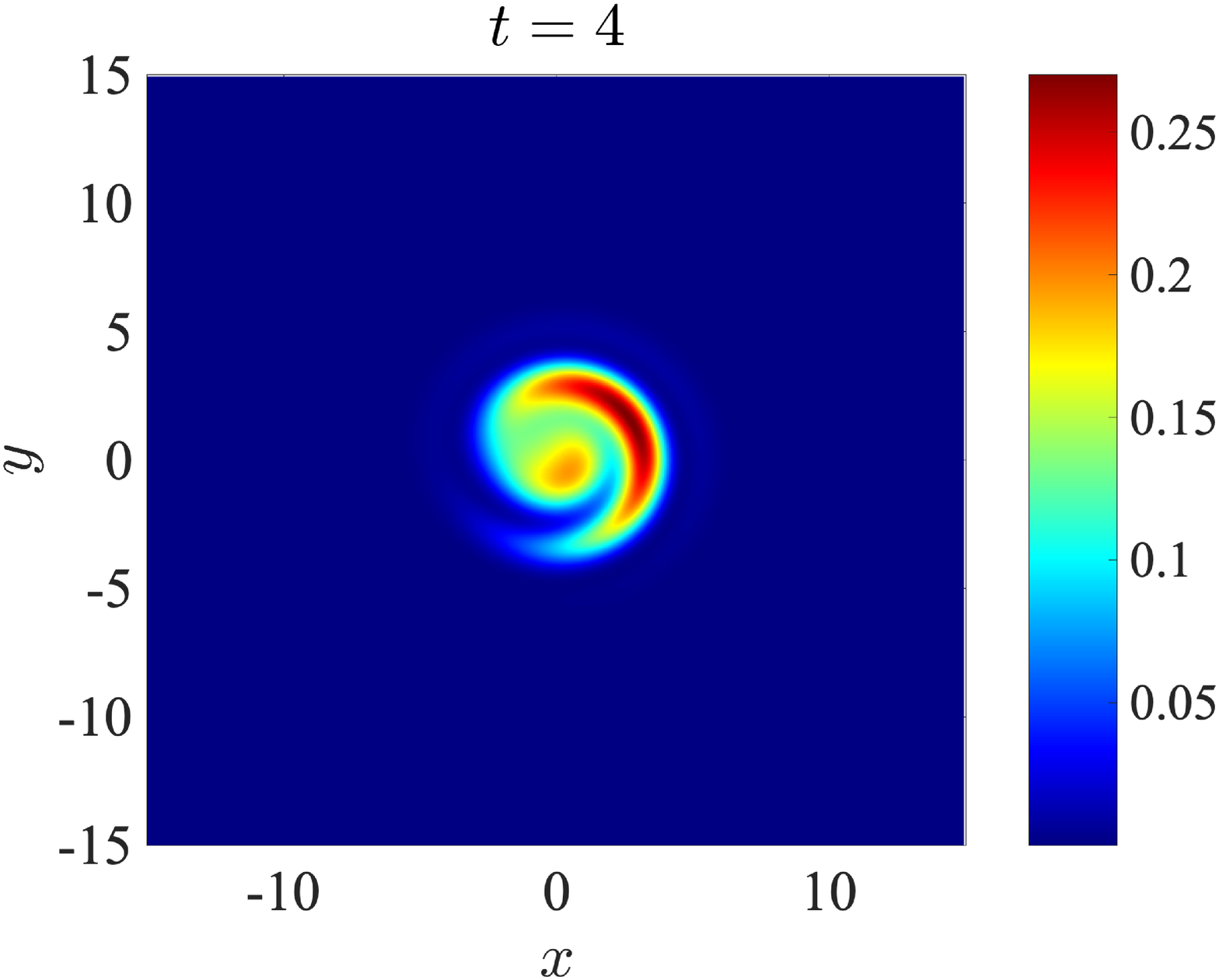}}
\end{minipage}
\begin{minipage}{0.5\textwidth}
\centerline{\includegraphics[width=8cm,height=6cm]{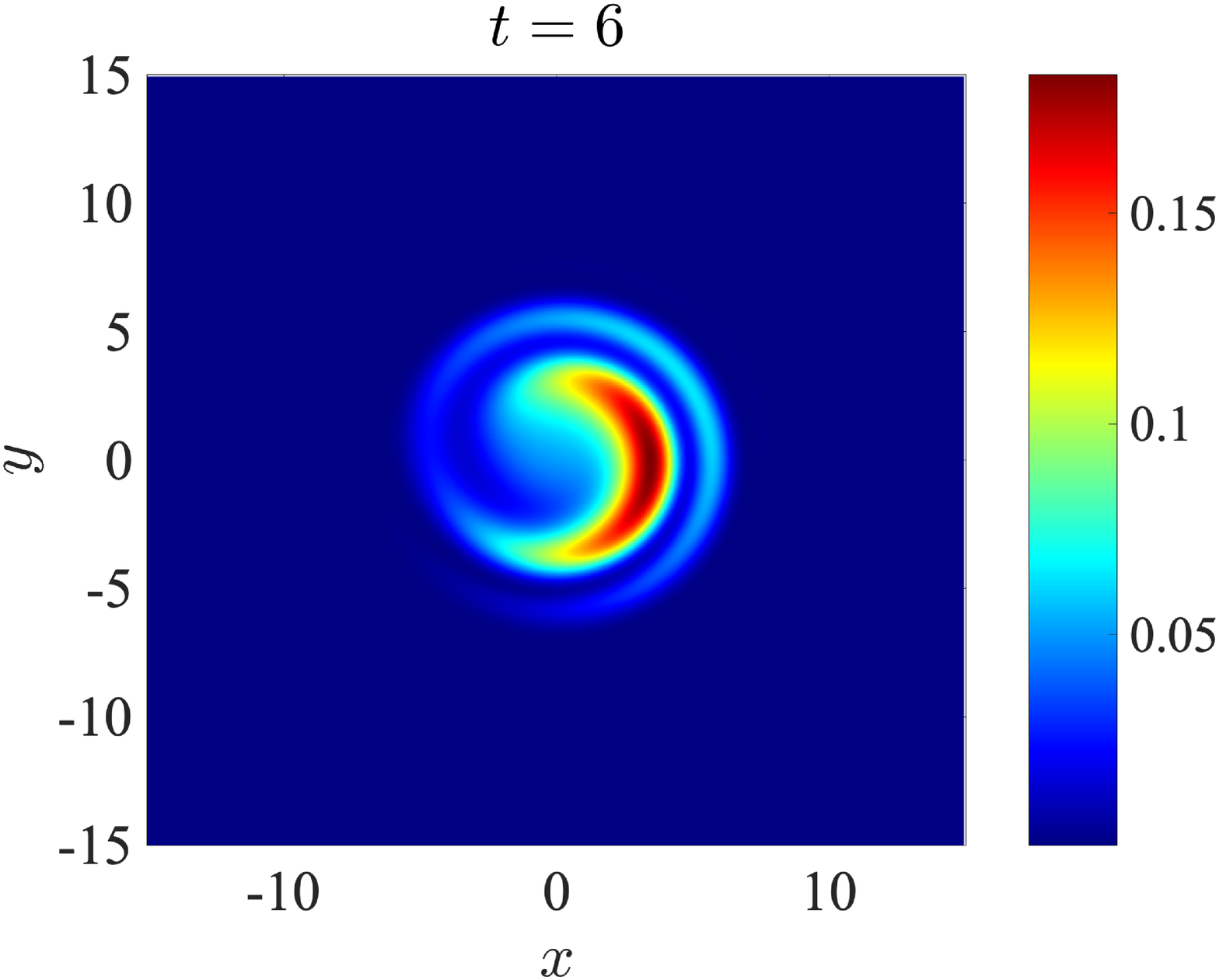}}
\end{minipage}
\begin{minipage}{0.5\textwidth}
\centerline{\includegraphics[width=8cm,height=6cm]{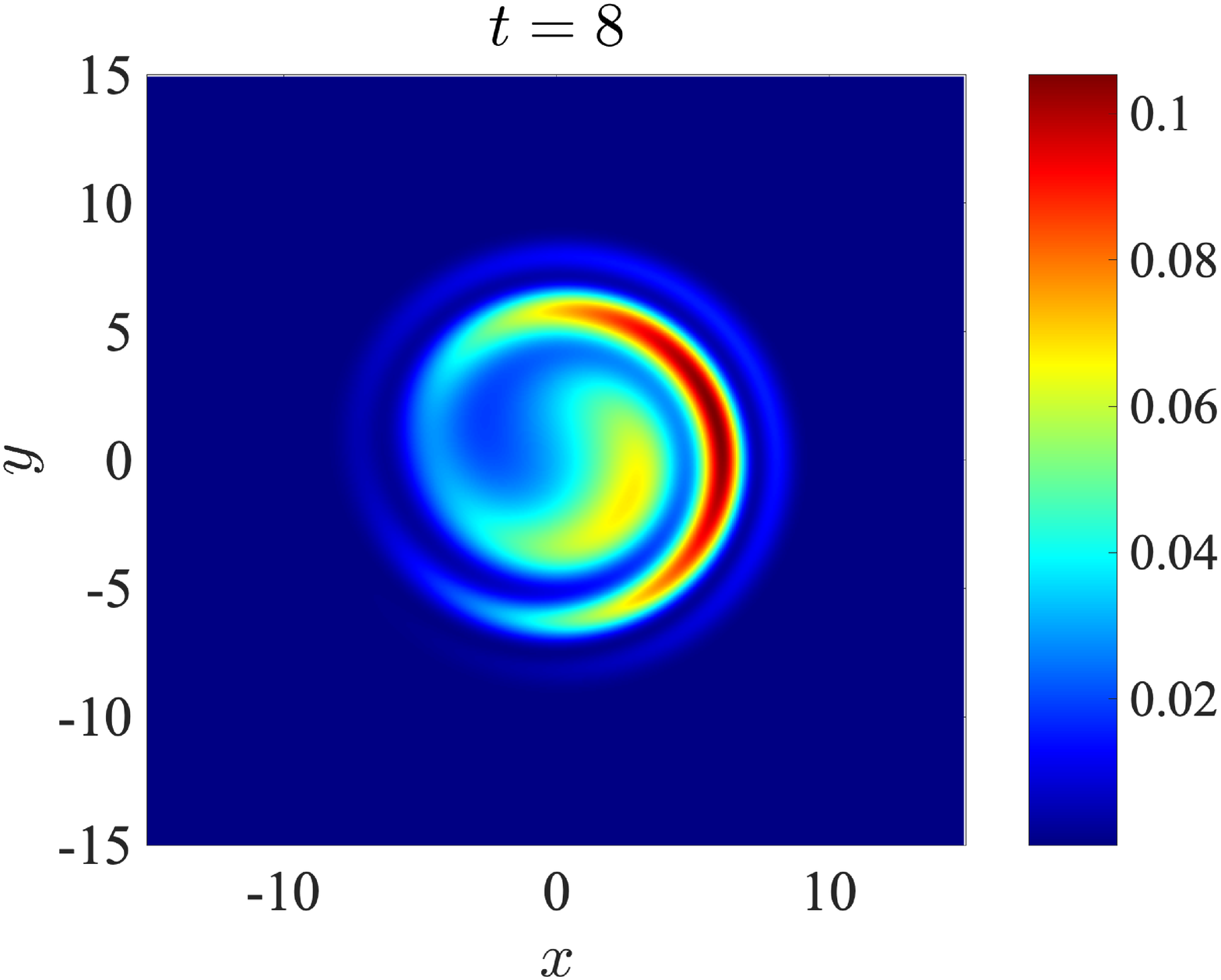}}
\end{minipage}
\begin{minipage}{0.5\textwidth}
\centerline{\includegraphics[width=8cm,height=6cm]{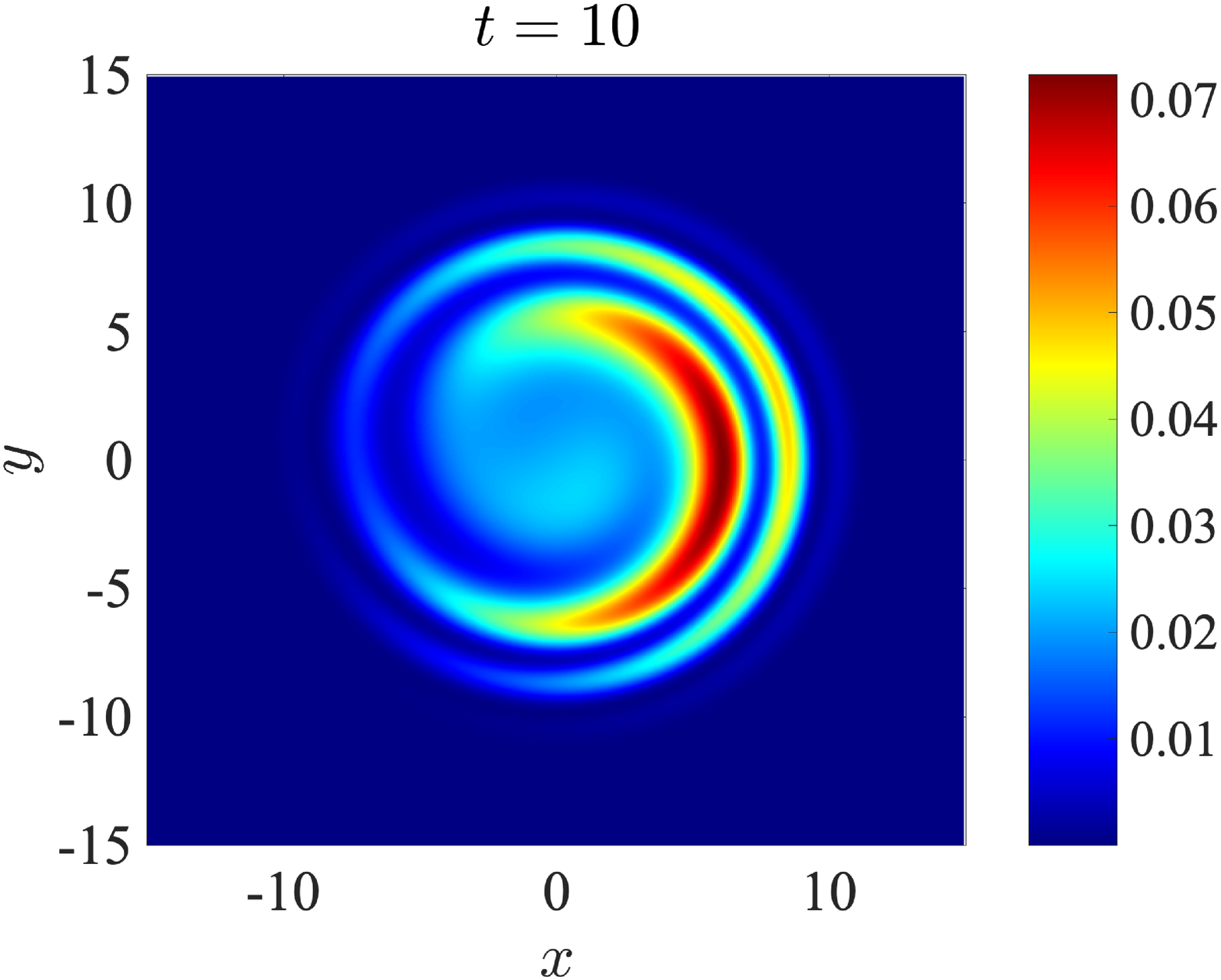}}
\end{minipage}
\caption{Dynamics of the density $\rho(t, {\bf x}) = |\phi_1(t, {\bf x})|^2 + |\phi_2(t, {\bf x})|^2 $ of the Dirac equation \eqref{eq:Dirac_21} in 2D with a honeycomb lattice potential when $\varepsilon = 0.001$ up to $t = 10$.}
\label{fig:eps0_rho}
\end{figure}

Figure \ref{fig:eps1_rho} and Figure \ref{fig:eps0_rho} depict the evolution from $t = 0$ to $t = 10$  of the density $\rho(t, {\bf x}) = |\phi_1(t, {\bf x})|^2 + |\phi_2(t, {\bf x})|^2$ of the Dirac equation \eqref{eq:Dirac_21} in 2D with honeycomb lattice potential when $\varepsilon = 1$ and $\varepsilon = 0.001$, respectively. From these two figures, we find out that the dynamics of the density depend on the small parameter $\varepsilon$ and the Zitterbewegung is observed forming a spiral, which is consistent with the experiment.  In addition, the EWI-FP/sEWI-FP method is able to capture the dynamics of the Dirac equation \eqref{eq:Dirac_21} in 2D accurately and efficiently.

\section{Conclusions}
The exponential wave integrator (EWI) methods combined with the Fourier spectral discretization in space were rigorously carried out and analyzed for the long-time dynamics of the Dirac equation with small potentials. The EWI-FP and sEWI-FP methods are explicit and the error bounds are at $O(h^{m_0} + \tau^2)$ up to the time at $O(1/\varepsilon)$. The $\varepsilon$-scalability of the EWI-FP and sEWI-FP methods up to the time at $O(1/\varepsilon)$ should be taken as $h = O(1)$ and $\tau = O(1)$, which is better than the finite difference time discretization in the long-time regime. The numerical results confirm our error estimates in the long-time regime and compare the long-time behaviors of the EWI-FP and sEWI-FP methods, which indicate that the symmetric scheme performs much better for the long-time simulations of the Dirac equation. Finally, we studied the dynamics of the Dirac equation in 2D with a honeycomb lattice potential and the Zitterbewegung was observed for the density with different $\varepsilon$.

\section*{Acknowledgements}

The authors would like to thank Professor Weizhu Bao for his valuable suggestions and comments. This work was partially supported by the Ministry of Education of Singapore grant R-146-000-290-114. Part of the work was done when the authors were visiting the Institute for Mathematical Sciences at the National University of Singapore in 2020. 



\begin{thebibliography}{}
%
     
    	 \bibitem{AH}
		{E. Ackad, M. Horbatsch},
	 	Numerical solution of the Dirac equation by a mapped Fourier grid method, J. Phys. A: Math. General 38 (2005) 3157--3171.      

    	 \bibitem{AL}
		{X. Antoine, E. Lorin},
	 	Computational performance of simple and efficient sequential and parallel Dirac equation solvers, Comput. Phys. Commun. 220 (2017) 150--172.   

   		\bibitem{BC}
		{W. Bao, Y. Cai},
	     Uniform and optimal error estimates of an exponential wave integrator sine pseudospectral method for the nonlinear Schr\"odinger equation with wave operator, SIAM J. Numer. Anal. 52 (2014) 1103--1127.          
	     
   		\bibitem{BCJT}
		{W. Bao, Y. Cai, X. Jia, Q. Tang},
	     Numerical methods and comparison for the Dirac equation in the nonrelativistic limit regime, J. Sci. Comput. 71 (2017) 1094--1134.
	
        \bibitem{BCJY}
        {W. Bao, Y. Cai, X. Jia, J. Yin},
        Error estimates of numerical methods for the nonlinear Dirac equation in the nonrelativistic limit regime, Sci. China Math. 59 (2016) 1461--1494.

        
         \bibitem{BCY}
        {W. Bao, Y. Cai, J. Yin},
          Uniform error bounds of time-splitting methods for the nonlinear Dirac equation in the nonrelativistic limit regime, SIAM J. Numer. Anal. 59 (2021) 1040--1066.    
        
         \bibitem{BSG}
        {J. W. Braun, Q. Su, R. Grobe},
          Numerical approach to solve the time-dependent Dirac equation, Phys. Rev. A 59 (1999) 604--612.
         
          \bibitem{BHM}
	     {D. Brinkman, C. Heitzinger, P. A. Markowich},
	      A convergent 2D finite-difference scheme for the Dirac-Poisson system and the simulation of graphene, J. Comput. Phys. 257 (2014) 318--332.
	      
          \bibitem{CW}
	     {Y. Cai, Y. Wang},
	      Uniformly accurate nested Picard iterative integrators for the Dirac equation in the nonrelativistic limit, SIAM J. Numer. Anal. 57 (2019) 1602--1624.	    

          \bibitem{CCO}
	     {E. Celledoni, D. Cohen, B. Owren},
	      Symmetric exponential integrators with an application to the cubic Schr\"odinger equation, Found. Comp. Math. 8 (2008) 303--317.	      	      
	      	      
          \bibitem{CC}
	     {R. J. Cirincione, P. R. Chernoff},
	      Dirac and Klein-Gordon equations: convergence of solutions in the nonrelativistic limit, Commun. Math. Phys. 79 (1981) 33--46.	      

          \bibitem{CG}
	     {D. Cohen, L. Gauckler},
	     One-stage exponential integrators for nonlinear Schr\"odinger equations over long times, BIT 52  (2011) 877--903.	      
         
         \bibitem{Das1}
        {A. Das},
        General solutions of {M}axwell-{D}irac equations in 1 + 1-dimensional space-time and spatially confined solution, J. Math. Phys. 34 (1993) 3986--3999.
        
         \bibitem{Das2}
        {A. Das, D. Kay},
        A class of exact plane wave solutions of the Maxwell-Dirac equations, J. Math. Phys. 30 (1989) 2280--2284.
        
         \bibitem{Dirac1}
        {P. A. M. Dirac},
        The quantum theory of the electron, Proc. R. Soc. Lond. A 117 (1928) 610--624.
        
         \bibitem{Dirac2}
        {P. A. M. Dirac},
        {P}rinciples of {Q}uantum {M}echanics, Oxford University Press, London, 1958.

         \bibitem{DF1}
        {G. Dujardin, E. Faou},
        Long time behavior of splitting methods applied to the linear Schr\"odinger equation, C. R. Acad. Sci. Paris 344 (2007) 89--92.        

        
         \bibitem{DF2}
        {G. Dujardin, E. Faou},
        Normal form and long time analysis of splitting schemes for the linear Schr\"odinger equation with small potential, Numer. Math. 108 (2007) 223--262.       
        
         \bibitem{ES}
        {M. Esteban, E. S\'er\'e},
        Existence and multiplicity of solutions for linear and nonlinear Dirac problems, Partial Differ. Equ. Appl. 12 (1997) 107--112.
                
         
         
         \bibitem{FY}
        {Y. Feng, J. Yin},
        Spatial resolution of different discretizations over long-time for the Dirac equation with small potentials.        
        
        
         \bibitem{FLB}
        {F. Fillion-Gourdeau, E. Lorin, A. D. Bandrauk},
        Resonantly enhanced pair production in a simple diatomic model, Phys. Rev. Lett. 110 (2013) 013002.
        
         \bibitem{FLB1}
        {F. Fillion-Gourdeau, E. Lorin, A. D. Bandrauk},
        Numerical solution of the time-dependent Dirac equation in coordinate space without fermion-doubling, Comput. Phys. Commun. 183 (2012) 1403--1415.
        
 	    \bibitem{GL}
	    {L. Gauckler, C. Lubich},
	    Splitting integrators for nonlinear Schr\"odinger equations over long times,
	    Found. Comput. Math., 10 (2010) 275--302.
 	
		\bibitem{Gau}
		{W. Gautschi},
		Numerical integration of ordinary differential equations based on trigonometric polynomials, Numer. Math. 3 (1961) 381--397.

	
         \bibitem{GGT}
        {F. Gesztesy, H. Grosse, B. Thaller},
        A rigorous approach to relativistic corrections of bound state energies for spin-1/2 particles, Ann. Inst. Henri Poincar\'e Phys. Theor. 40 (1984) 159--174.


         \bibitem{Gosse}
        {L. Gosse},
        A well-balanced and asymptotic-preserving scheme for the one-dimensional linear Dirac equation, BIT 55 (2015) 433--458.
        
  		 \bibitem{Gr1}
		{V. Grimm},
		A note on the Gautschi-type method for oscillatory second-order differential equations, Numer. Math. 102 (2005) 61--66.
	
		\bibitem{Gr2}
		{V. Grimm},
		On error bounds for the Gautschi-type exponential integrator applied to oscillatory second-order differential equations, Numer. Math. 100 (2005) 71--89.
	
		\bibitem{GH}
		{V. Grimm, M. Hochbruck},
		Error analysis of exponential integrators for oscillatory second-order differential equations, J. Phys. A, 39 (2006), 5495--5507.
               
         \bibitem{Gross}
        {L. Gross},
        The Cauchy problem for the coupled Maxwell and Dirac equations, Commun. Pure Appl. Math. 19 (1966) 1--15.
        
         \bibitem{GSX}
        {B.-Y. Guo, J. Shen, C.-L. Xu},
        Spectral and pseudospectral approximations using Hermite functions: Application to the Dirac equation, Adv. Comput. Math. 19 (2003) 35--55.

         
         \bibitem{HLW}
	     {E. Hairer, C. Lubich, G. Wanner},
	     Geometric Numerical Integration,
	     Springer, New York, 2002.
	                 
         \bibitem{HL}
        {M. Hochbruck, C. Lubich},
        Exponential integrators for quantum-classical molecular dynamics, BIT 39 (1999) 620--645.       
 
          \bibitem{HLS}
        {M. Hochbruck, C. Lubich, H. Selhofer},
        Exponential integrators for large systems of differential equations, SIAM J. Sci. Comput. 19 (1998) 1552--1574.       

               
         \bibitem{HO}
        {M. Hochbruck, A. Ostermann},
        Exponential integrators, Acta Numer. 19 (2000) 209--286.       
              
         
         \bibitem{MY}
        {Y. Ma, J. Yin},
         Error bounds of the finite difference time domain methods for the Dirac equation in the semiclassical regime,  J. Sci. Comput 81 (2019) 1801--1822.
      

         
         \bibitem{ST}
	     {J. Shen, T. Tang},
	     Spectral and High-Order Methods with Applications,
	     Science Press, Beijing, 2006.
         
         \bibitem{Smith}
        {G. D. Smith},
        Numerical Solution of Partial Differential Equations: Finite Difference Methods, Clarendon Press, Oxford (1985).
        
        
        \bibitem{WHJY}
        {H. Wu, Z. Huang, S. Jin, D. Yin},
         Gaussian beam methods for the Dirac equation in the semi-classical regime, Commun. Math. Sci. 10 (2012) 1301-1315.
         

\end{thebibliography}
\end{document}